\documentclass[final]{siamltex}
\usepackage{graphicx}
\usepackage{caption}
\usepackage{subfigure}
\usepackage{amsmath}
\usepackage{amsfonts}
\usepackage{amssymb}
\usepackage{enumerate}
\usepackage[ruled,vlined]{algorithm2e}
\usepackage{url}
\usepackage{pgf,pgfarrows,pgfnodes,pgfautomata,pgfheaps,pgfshade}
\usepackage{booktabs,multirow,lscape}

\newcommand{\R}{\mathbb{R}}

\begin{document}
\title{Using first-order information in Direct Multisearch for
multiobjective optimization}
\author{R. Andreani
\thanks{Department of Applied Mathematics, IMECC-UNICAMP, University  of Campinas, Rua S\'ergio Buarque de Holanda, 651 Cidade Universit\'aria ``Zeferino Vaz", Distrito Bar\~ao Geraldo,
13083-859 Campinas SP, Brazil (\texttt{andreani@unicamp.br}). This
author was financially supported by FAPESP (Projects 2013/05475-7
and 2017/18308-2)  and CNPq (Project 301888/2017-5).} \and A. L.
Cust\'odio
\thanks{Department of Mathematics, FCT-UNL-CMA, Campus de Caparica, 2829-516 Caparica, Portugal ({\tt alcustodio@fct.unl.pt}). Support for this author was
provided by national funds through FCT -- Funda\c c\~ao para a Ci\^encia e a Tecnologia I. P., under the scope of projects
PTDC/MAT-APL/28400/2017 and UIDB/MAT/00297/2020.}
\and M. Raydan
\thanks{Centro de Matem\'atica e Aplica\c{c}\~oes (CMA), FCT, UNL, 2829-516 Caparica, Portugal ({\tt m.raydan@fct.unl.pt}). This author was
financially supported by  the Funda\c c\~ao para a Ci\^encia e a Tecnologia (Portuguese Foundation for Science and Technology) through the projects
  UIDB/MAT/00297/2020 (Centro de Matem\'atica e Aplica\c{c}\~oes) and  CEECIND/02211/2017. } }

\date{January 17, 2021}
\maketitle
\begin{abstract}
Derivatives are an important tool for single-objective
optimization. In fact, it is commonly accepted that
derivative-based methods present a better performance than
derivative-free optimization approaches. In this work, we will
show that the same does not apply to multiobjective
derivative-based optimization, when the goal is to compute an
approximation to the complete Pareto front of a given problem. The
competitiveness of Direct MultiSearch (DMS), a robust and
efficient derivative-free optimization algorithm, will be stated
for derivative-based multiobjective optimization problems. We will
then assess the potential enrichment of adding first-order
information to the DMS framework. Derivatives will be used to
prune the positive spanning sets considered at the poll step of
the algorithm, highlighting the role that ascent directions, that
conform to the geometry of the nearby feasible region, can have.
Both variants of DMS show to be competitive against a state-of-art
derivative-based algorithm. Moreover, for reasonable small budgets
of function evaluations, the new variant is not only competitive
with the derivative-based solver but also with the original
implementation of DMS.
\\ [2mm]
[2mm] {\bf Keywords:}  Multiobjective optimization, Pareto front
computation, Derivative-based methods, Derivative-free
optimization, Direct MultiSearch.

\end{abstract}

\section{Introduction}

Multiobjective optimization (MOO) problems appear frequen-tly in
engineering and scientific applications, in such diverse areas
such as civil engineering, environment, medicine or aerospace
engineering~\cite{HAfshari_et_al_2019,FMaglia_et_al_2018,PPotrebko_et_al_2017,ARavanbakhsh_SFranchini_2012},
just to cite a few. The major feature of a MOO problem is the
presence of finitely many components in the objective function,
associated to conflicting objectives, that have to be
simultaneously optimized. Hardly a single point will optimize all
of them at once, hence a nonstandard notion of optimality is
required. The fundamental optimality concept is that of Pareto
optimal point, which is a point such that no improvement in all
the components of the objective function can be achieved by moving
to another feasible point. The image set of all Pareto optimal
points (also called the Pareto front) is usually a continuum that
may have disjoint components. In general, for a problem with $p>1$
objectives, the Pareto front is a manifold of dimension $p-1$. For
example, if $p=2$ the Pareto front will be a curve (or a set of
curve segments), which provides in a compact way all the
information required for a user to choose an appropriate Pareto
optimal point as a compromise solution between the usually
conflicting components of the objective function. Like in
classical single objective optimization, finding global Pareto
optimal points is difficult, unless additional information is
available about the objective function. Thus, MOO algorithms
typically try to find local Pareto optimal points for the
problems, meaning that the definition of Pareto optimality is
satisfied in a neighborhood of the current point.

There are several classes of MOO algorithms, depending not only on
the level of smoothness of the objective function but also on the
time when the user establishes an order preference for the
different components of the objective
function~\cite{RTMarler_JSArora_2004}. In this work, we will focus
on methods with \emph{a posteriori articulation of preferences},
which attempt to capture the whole Pareto front of the problem,
never establishing preferences among the several components of the
objective function. Evolutionary algorithms, or other similar
heuristics, belong to this class. However, these algorithms lack a
well-established convergence analysis and are usually slow in
converging to the Pareto front of the problem, requiring a large
number of iterations and function
evaluations~\cite{MTMEmmerich_AHDeutz_2018}. When derivatives of
the different components of the objective function are available,
typical approaches are based on generating a single sequence of
iterates that converges to a point with corresponding image lying
on the Pareto front (one at a time); see,
e.g.,~\cite{GACarrizo_et_al_2016, NSFazzio_MLSchuverdt_2019,
JFliege_et_al_2009, JFliege_BFSvaiter_2000,
LMGranaDrummond_ANIusem_2004}. Multistart
approaches~\cite{EMiglierina_et_al_2008} or scalarization
techniques~\cite{GEichfelder_2008}, can help in finding
approximations to the complete Pareto front of a given MOO
problem. Although, the first can be computational expensive and
the latter generally fails in detecting nonconvex parts of
it~\cite{IDas_JEDennis_1998}.

Recently, a novel approach has been developed to approximate the
entire Pareto front using first and second order
information~\cite{JFliege_AIFVaz_2016}. The so-called MOSQP method
keeps a list of nondominated points, which approximates the Pareto
front of the MOO problem, that is improved both for spread along
the Pareto front and optimality by solving single-objective
constrained optimization problems derived as SQP problems.
In~\cite{JFliege_AIFVaz_2016}, numerical results are reported
indicating the superiority of the MOSQP algorithm when compared to
other state-of-the-art multiobjective solvers.

In derivative-free optimization, Direct MultiSearch
(DMS)~\cite{ALCustodio_et_al_2011} is also able to compute
approximations to the complete Pareto front of a given MOO
problem. This is a well-established algorithm, with theoretical
results regarding convergence, and consistently used with good
results both for benchmark of new
solvers~\cite{GCocchi_et_al_2018,GLiuzzi_SLucidi_FRinaldi_2016} or
in real applications~\cite{RBrito_et_al_2019,DHirpa_et_al_2016}.

The purpose of the current work is twofold. Our first objective is
to compare the performance of the derivative-free DMS method and
the derivative-based MOSQP algorithm. In single objective
optimization, it is common to say that if derivatives are
available, or can be obtained at a reasonable cost (e.g. using
finite-differences) then derivative-based optimization is
preferable to derivative-free optimization methods~\cite[p.
6]{CAudet_WHare_2017}. We will provide numerical results on a
large set of benchmark MOO problems, that allow to assess the
numerical performance of derivative-based and derivative-free
optimization solvers, when computing approximations to the
complete Pareto fronts of derivative-based MOO problems. Our
second objective is to asses the potential enrichment of adding
first-order information, when derivatives are available, to the
DMS framework. We will describe and analyze several different
combined techniques that maintain the search/poll paradigm of DMS,
while adding in a convenient way gradient information to the poll
step. Again, the value of the proposed strategies will be assessed
through numerical experiments.

The remaining of this document is organized as follows. In
Section~\ref{Sec:DMS} we present the MOO problem and briefly
revise the derivative-free optimization method DMS, since it will
be later modified to incorporate first-order information.
Section~\ref{Sec:DMSvsMOSQP} is devoted to a full numerical
comparison between DMS and MOSQP methods.
Section~\ref{Sec:Pruning} details the use of first-order
information to eliminate directions in the poll step of DMS and
assesses the corresponding numerical performance. In
Section~\ref{Sec:Ascent} the usefulness of ascent directions is
motivated by illustrating their performance on one properly chosen
biobjective problem. Results are then reported in the complete
test set. Finally, in Section~\ref{Sec:Conclusions}
 we present some concluding remarks.

\section{DMS at a glance}\label{Sec:DMS}

We consider multiobjective minimization problems of the form

\begin{equation} \label{probm}
\min_{x\in\Omega} F(x)\equiv (f_1(x),\dots,f_p(x))^{\top},
\end{equation}

\noindent where $p\geq 2$, $\Omega\subseteq \R^{n}$ represents the
feasible region, typically defined as a box
$\Omega=\{x\in\mathbb{R}^n:l \leq x \leq u\}$, and for each $i$
($1\leq i\leq p$) $f_i: \Omega \rightarrow \R \cup \{+\infty\}$
denotes a component of the objective function, which we assume to
be strictly differentiable in $\Omega$ (continuity of the partial
derivatives is not required).

The Direct MultiSearch (DMS) method was originally proposed
in~\cite{ALCustodio_et_al_2011}, generalizing directional direct
search to multiobjective derivative-free optimization. For a
review on single objective derivative-free optimization methods we
recommend~\cite{CAudet_WHare_2017, ARConn_et_al_2009}. The
algorithm has also been successfully extended to global
multiobjective derivative-free
optimization~\cite{ALCustodio_JFAMadeira_2018}, by coupling it
with a multistart initialization technique, where not all the
initialized searches are conducted until the end.

Being a directional direct search method, each iteration of DMS
conforms to the search/poll paradigm. The search step is optional,
since the convergence results derive from the poll step of the
algorithm. In fact, in the original presentation of the
method~\cite{ALCustodio_et_al_2011}, it was left empty and this
will be the approach followed in the present work. Recently, the
minimization of quadratic polynomial models, which have always
played a key role in derivative-free methods for single objective
optimization, was used for successfully defining a search step for
DMS~\cite{CPBras_ALCustodio_2020}. First-order information can
surely be used to define appropriate search steps, following the
strategies proposed in~\cite{CPBras_ALCustodio_2020}, but that
will not be the subject of the present work, which will focus on
the poll step.

We present a simplified description of the DMS framework, where
only the poll step is considered, and where the globalization
strategy is based on the use of integer lattices, meaning that all
the points generated by the algorithm lie on an implicit mesh. For
a more general description, we refer to the original
work~\cite{ALCustodio_et_al_2011}.

The algorithm initializes with a list of feasible, nondominated
points (possibly just one) and corresponding stepsize parameters.
Making use of the strict partial order induced by the cone
$\mathbb{R}_{+}^p$, we say that point $x$ dominates point $y$ when
$F(x)\prec_F F(y)$, i.e., when $F(y)-F(x) \in
\mathbb{R}_{+}^p\setminus\{0\}$. If $x$ does not dominate $y$ and
$y$ does not dominate $x$, $x$ and $y$ are said to be
nondominated. The list, representing the current approximation to
the Pareto front of the MOO problem, will be updated at every
iteration by generating new feasible points which are compared
with the points already stored in it, only keeping the
nondominated ones.

At each iteration, a feasible nondominate point stored in the list
and the associated stepsize parameter, will be selected. Different
strategies can be considered in the selection of this poll center.
Currently it is based on a spread
metric~\cite{ALCustodio_et_al_2011}, in an attempt of reducing the
gaps between consecutive points lying in the current approximation
to the Pareto front of the problem.

The poll step of the algorithm consists on a local search around
the selected poll center, by testing a set of directions with an
adequate geometry, scaled by the corresponding stepsize parameter.
Typically, positive spanning sets are
considered~\cite{CDavis_1954}, that should conform to the geometry
of the nearby active constraints of the current poll
center~\cite{TGKolda_RMLewis_VTorczon_2003}.

For convergence purposes, the poll step can be performed either in
a complete or an opportunistic way. In the latter, the polling
procedure is stopped once a new feasible nondominated point is
found. The complete approach tests all the poll directions, only
adding to the list the new feasible nondominated points found (and
removing from the list all the dominated ones). We will follow
this last approach, which is the one corresponding to the original
algorithmic implementation of DMS~\cite{ALCustodio_et_al_2011}, in
an attempt of maximizing the number of feasible nondominated
points generated at each iteration.

The final step of each iteration is the update of the stepsize
parameter, which is increased or kept constant for successful
iterations and decreased for unsuccessful ones. An iteration is
said to be successful if the list changes, meaning that at least
one new feasible nondominated point was found. Unsuccessful
iterations keep the list unchanged.

A simplified description of the DMS framework is provided in
Algorithm~\ref{alg:dms}. For a complete description see~\cite{ALCustodio_et_al_2011}.\\

\begin{algorithm}
\begin{rm}
\begin{description}
\item[Initialization] \ \\
Choose a set of nondominated points $\{x_{ini}^i\in\Omega,\, i\in
I\}$ with $f_j(x_{ini}^i)<+\infty,\forall j
\in\{1,\ldots,p\},\forall\, i \in I$, $\alpha_{ini}^i>0,\, i\in I$
initial stepsizes, $0 < \beta_1 \leq \beta_2 < 1$ the coefficients
for stepsize contraction and $\gamma \geq 1$ the coefficient for
stepsize expansion. Let $\mathcal{D}$ be a set of positive
spanning sets. Initialize the list of feasible nondominated points
and corresponding stepsize parameters $L_0=\{
(x_{ini}^i;\alpha_{ini}^i ), i\in I \}$. \vspace{1ex}
\item[For $k=0,1,2,\ldots$] \ \\
\begin{enumerate}
\item {\bf Selection of an iterate point:} Order the list $L_k$
according to some criteria and select the first item $(x;\alpha)
\in L_k$ as the current iterate and stepsize parameter (thus
setting $(x_k;\alpha_k)=(x;\alpha)$). \item {\bf Poll step:}
Choose a positive spanning set~$D_k$ from the set $\mathcal{D}$.
Evaluate $F$ at the feasible poll points belonging to $\{
x_k+\alpha_k d: \, d \in D_k \}$. Compute $L_{trial}$ by removing
all dominated points from $L_{k} \cup \{(x_k+\alpha_k d;\alpha_k):
d\in D_k \wedge x_k+\alpha_k d \in \Omega\}$. If $L_{trial} \neq
L_k$ declare the iteration (and the poll step) successful and set
$L_{k+1}=L_{trial}$. Otherwise, declare the iteration (and the
poll step) unsuccessful and set $L_{k+1} = L_k$. \vspace{1ex}
\item {\bf Stepsize parameter update:} If the iteration was
successful then maintain or increase the corresponding stepsize
parameters, by considering $\alpha_{k,new}$ $\in [\alpha_k, \gamma
\alpha_k]$ and replacing all the new points $(x_k+\alpha_k
d;\alpha_k)$ in $L_{k+1}$ by $(x_k + \alpha_k d;\alpha_{k,new})$.
Replace also $(x_k;\alpha_k)$, if in $L_{k+1}$, by
$(x_k;\alpha_{k,new})$.\\
Otherwise, decrease the stepsize parameter, by choosing
$\alpha_{k,new} \in [\beta_1 \alpha_k, \beta_2 \alpha_k]$, and
replace the poll pair $(x_k;\alpha_k)$ in $L_{k+1}$ by
$(x_k;\alpha_{k,new})$. \vspace{1ex}
\end{enumerate}
\end{description}
\end{rm}
\caption{A simplified description of Direct MultiSearch (DMS).}
\label{alg:dms}
\end{algorithm}

The convergence of DMS has been established
in~\cite{ALCustodio_et_al_2011}, closely following the arguments
used in the analysis of single objective directional direct search
methods. After stating the existence of a subsequence of stepsize
parameters converging to zero, this property is used for
establishing Pareto-Clarke-KKT criticality. The result is
formalized in Theorem~\ref{th:dense} for limit points of
convergent refining subsequences.

\begin{definition} A subsequence $\{x_k\}_{k\in K}$ of iterates
corresponding to unsuccessful poll steps is said to be a refining
subsequence if $\{\alpha_k\}_{k\in K}$ converges to zero.
\end{definition}

The concept of refining direction is associated with convergent
refining subsequences and is formalized in
Definition~\ref{def:refining directions}.

\begin{definition}\label{def:refining directions}
Let $x_*$ be the limit point of a convergent refining subsequence
$\{x_k\}_{k\in K}$. If the limit $\lim_{k\in K'} d_k/\|d_k\|$
exists, where $K'\subseteq K$ and $d_k\in D_k$, and if
$x_k+\alpha_k d_k \in \Omega$, for sufficiently large $k \in K'$,
then this limit is said to be a refining direction for~$x_*$.
\end{definition}

Assuming the density of the set of refining directions in the
Clarke tangent cone to $\Omega$ computed at limit points of
refining subsequences~\cite{FHClarke_1990}, the convergence of DMS
is established.

\begin{theorem}(see~\cite{ALCustodio_et_al_2011})\label{th:dense}
Consider a refining subsequence $\{x_k\}_{k\in K}$ converging to
$x_* \in \Omega$. Assume that $F$ is strictly differentiable at
$x_*$ and that the interior of the tangent cone to $\Omega$ at
$x_*$ is nonempty. If the set of refining directions for $x_*$ is
dense in the Clarke tangent cone to $\Omega$ at $x_*$, then $x_*$
is a Pareto-Clarke-KKT critical point, i.e,
$$\forall d\in T_\Omega^{Cl}(x_*),\exists
j(d)\in\{1,2,\ldots,p\}:\nabla f_{j(d)}(x_*)^\top d \geq 0.$$
\end{theorem}

Recently, worst-case complexity bounds were provided for DMS, but
considering a globalization strategy that requires sufficient
decrease for accepting new points~\cite{ALCustodio_et_al_2021}.
For a particular algorithmic instance, which considers a stricter
criterion for accepting new nondominated points, DMS presents a
worst-case complexity bound of $\mathcal{O}(\epsilon^{-2})$.
similar to the one of steepest descent.

\section{Comparing DMS and MOSQP}\label{Sec:DMSvsMOSQP}

Derivatives are a keystone for optimization. As previously
mentioned, in single objective optimization, when in presence of
smooth functions, derivative-based methods are preferable to
derivative-free optimization algorithms, even if one has to spend
some time and effort to obtain good quality derivatives
(see~\cite[p. 7]{ARConn_et_al_2009} or~\cite[p.
6]{CAudet_WHare_2017}). In this section, we will assess the
situation for multiobjective derivative-based optimization, when
the goal is to compute approximations to complete Pareto fronts.

For that, DMS algorithm~\cite{ALCustodio_et_al_2011} will be
numerically tested against MOSQP~\cite{JFliege_AIFVaz_2016}. The
latter is a recent solver proposed for multiobjective
derivative-based optimization which uses a SQP approach. MOSQP
keeps a list of nondominated points that is improved both for
spread along the Pareto front and optimality by solving
single-objective constrained optimization problems. Thus, MOSQP is
able to generate approximations to complete Pareto fronts, an
advantage over classical derivative-based multiobjective
optimization solvers, which compute a single point. At the time of
the release, extensive numerical results were provided for MOSQP,
including a comparison with a classical scalarization approach for
biobjective problems~\cite{JFliege_AIFVaz_2016}. The good results
obtained, allowed the authors to conclude that MOSQP  should be
\emph{``the preferred solution framework for multiobjective
optimization problems when derivatives of objective and constraint
functions are available''}~\cite{JFliege_AIFVaz_2016}, which
justifies our algorithmic choice as baseline against DMS. Default
parameters were considered for both solvers, with exception to the
maximum number of function evaluations allowed, which was set to
$20\,000$. In some cases, a small budget of $500$ function
evaluations was additionally considered, to ensure that the
conclusions drawn were not dependent on the large number of
function evaluations allowed.

As test set, we considered the collection of 100 bound constrained
multiobjective optimization problems available at
\texttt{http://www.mat.uc.pt/dms}. This collection was previously
used to test DMS and MOSQP, at the time of their first
release~\cite{ALCustodio_et_al_2011,JFliege_AIFVaz_2016}. From
this collection, we selected a total of $54$ problems, for which
we were able to guarantee the existence of derivatives.
Table~\ref{tab:probs} reports the resulting test set, which
comprises problems with $2$ or $3$ components in the objective
function and a number of variables, $n$, between $1$ and $30$.

\begin{table}
\centering
\begin{tabular}{lrr|lrr|lrr|lrr}
\hline Problem&$n$&$p$&Problem&$n$&$p$&Problem&$n$&
$p$&Problem&$n$&$p$\\
\hline BK1&2&2&DTLZ4n2&2&2&lovison3&2&2&MOP7&2&3\\
CL1&4&2&DTLZ6&22&3&lovison4&2&2&SK1&1&2\\
Deb41&2&2&DTLZ6n2&2&2&lovison5&3&3&SK2&4&2\\
Deb513&2&2&ex005&2&2&lovison6&3&3&SP1&2&2\\
Deb521b&2&2&Far1&2&2&LRS1&2&2&SSFYY1&2&2\\
DG01&1&2&Fonseca&2&2&MHHM1&1&3&SSFYY2&1&2\\
DPAM1&10&2&IKK1&2&3&MHHM2&2&3&TKLY1&4&2\\
DTLZ1&7&3&IM1&2&2&MLF1&1&2&VFM1&2&3\\
DTLZ1n2&2&2&Jin1&2&2&MLF2&2&2&VU1&2&2\\
DTLZ2&12&3&Jin3&2&2&MOP1&1&2&VU2&2&2\\
DTLZ2n2&2&2&L2ZDT2&30&2&MOP2&4&2&ZDT2&30&2\\
DTLZ3&12&3&L3ZDT2&30&2&MOP3&2&2&ZLT1&10&3\\
DTLZ3n2&2&2&lovison1&2&2&MOP5&2&3&&&\\
DTLZ4&12&3&lovison2&2&2&MOP6&2&2&&&\\
\hline
\end{tabular}
\caption{The test set considered in the numerical experiments.
Here~$n$ represents the number of variables and~$p$ is the number
of components of the objective function.\label{tab:probs}}
\end{table}

For performance assessment, we considered typical metrics from the
multiobjective optimization literature, like is the case of purity
and spread metrics, as defined in~\cite{ALCustodio_et_al_2011},
and also the hypervolume
indicator~\cite{EZitzler_LThiele_1998,EZitzler_et_al_2003}. While
other metrics could have been
considered~\cite{CAudet_et_al_2020,SCheng_et_al_2012,NRiquelme_et_al_2015},
these are typical choices in recent
works~\cite{MATAnsary_GPanda_2021,CPBras_ALCustodio_2020,GCocchi_et_al_2018}.

In a simplified view, purity measures the percentage of nondominated points generated by a given solver. For problem
$\hat{p}\in \mathcal{P}$ and solver $s\in \mathcal{S}$, purity is defined by the ratio
\begin{equation*}
\bar{t}_{\hat{p},s} = \frac{| F_{\hat{p},s} \cap F_{\hat{p}} |}{ |
F_{\hat{p},s} |},
\end{equation*}
where $F_{\hat{p},s}$ denotes the approximation to the Pareto front computed for problem $\hat{p}\in \mathcal{P}$ by solver
 $s\in \mathcal{S}$ and $F_{\hat{p}}$ is a reference Pareto front for  problem $\hat{p}\in P$. This reference Pareto front is computed by
 joining the final approximations computed by any of the solvers  tested and removing from it all the dominated points. A value of
 purity near one indicates that the majority of the points generated by the corresponding solver is nondominated.
However, these could be concentrated in a single part of the true Pareto front. Spread metrics are required to have a fair
assessment of the solver's performance.

Since the goal is to build an approximation to the complete Pareto front of each problem, the computation of spread metrics initiates
with the computation of the so-called `extreme points' of the Pareto front (see~\cite{ALCustodio_et_al_2011}). The spread
$\Gamma$ measures the maximum gap between consecutive points lying in the approximated Pareto front. The metric $\Gamma_{\hat{p},s} > 0$ for
 problem $\hat{p}\in \mathcal{P}$ and solver  $s\in \mathcal{S}$ is given by
\begin{equation}\label{gamma_metric}
 \Gamma_{\hat{p},s} \; = \; \max_{j \in
\{1,\dots,p\}}\left(\max_{i\in
\{0,\dots,N\}}\{\delta_{i,j}\}\right),
\end{equation}
where $\delta_{i,j}=(f_{j}(x_{i+1})-f_{j}(x_i))$,
$x_1,x_2,\ldots,x_N$ represent the points generated by solver $s \in {\cal S}$ for problem $\hat{p}\in \mathcal{P}$, and $x_0,
x_{N+1}$ correspond to the `extreme points'. Implicitly, we are assuming that the objective function values have been sorted by
increasing order for each objective~$j\in \{1,\dots,p\}$.

The spread metric $\Delta$~\cite{KDeb_et_al_2002} measures the uniformity of the gaps across the approximation to the Pareto front:
\begin{equation}\label{delta_metric} \Delta_{\hat{p},s} \; = \;
\max_{j\in\{1,\dots,p\}}\left(\frac{\delta_{0,j}+\delta_{N,j}+\sum_{i=1}^{N-1}|\delta_{i,j}-
\bar{\delta}_j|}{\delta_{0,j}+\delta_{N,j}+(N-1)\bar{\delta}_j}\right),
\end{equation}
\noindent where $\bar{\delta}_j$, for $j=1,\ldots,p$, represents the average of the distances $\delta_{i,j}$, $i=1,\dots,N-1$.

The fourth metric considered is the hypervolume indicator~\cite{EZitzler_et_al_2003}, which measures the volume of
the portion of the objective function space that is dominated by the computed approximation to the Pareto front of the problem, and
upper bounded by a given reference point $U_{\hat{p}} \in \mathbb{R}^p$. This reference point should be dominated by all
points belonging to the approximations computed for the Pareto front of a given problem $\hat{p}\in \mathcal{P}$. Formally, it
can be defined as:

$$HV_{\hat{p},s} =  Vol\{y \in \mathbb{R}^p| \, y \le U_{\hat{p}} \wedge \exists x \in F_{\hat{p},s} : x \le y\} = Vol \left(\bigcup_{x \in F_{\hat{p},s}}
 [x, U_{\hat{p}}]\right),$$

\noindent where $Vol(.)$ denotes the Lebesgue measure of a
$p$-dimensional set of points and $[x, U_{\hat{p}}]$ denotes the
interval box with lower corner $x$ and upper corner $U_{\hat{p}}.$
The approach proposed in~\cite{CMFonseca_et_al_2006} was used for
its practical computation and the resulting hypervolume values
were scaled to the interval $[0,1]$, following the procedure
described in~\cite{CPBras_ALCustodio_2020}.

Performance profiles~\cite{EDDolan_JJMore_2002} will be depicted for each of the four metrics considered. Let $t_{\hat{p},s}$
denote the performance of solver~$s \in \mathcal{S}$ on problem~$\hat{p}\in \mathcal{P}$, assuming that lower values of
$t_{\hat{p},s}$ indicate a better performance. Each performance profile represents the curve
$$\rho_s(\tau)=\frac{1}{|\mathcal{P}|}|\{\hat{p}\in
\mathcal{P}:r_{\hat{p},s}\leq\tau\}|,$$ with
$r_{\hat{p},s}=t_{\hat{p},s}/\min\{t_{\hat{p},\bar{s}}:\bar{s}\in\mathcal{S}\}$.
In the case of purity and hypervolume metrics, larger values
indicate better performance. Thus, when computing performance
profiles for these two metrics, we set
$t_{\hat{p},s}=1/t_{\hat{p},s}$, as proposed
in~\cite{ALCustodio_et_al_2011}.

\begin{figure}[htbp]
\begin{center}
\includegraphics[width=6.4cm,height=4cm]{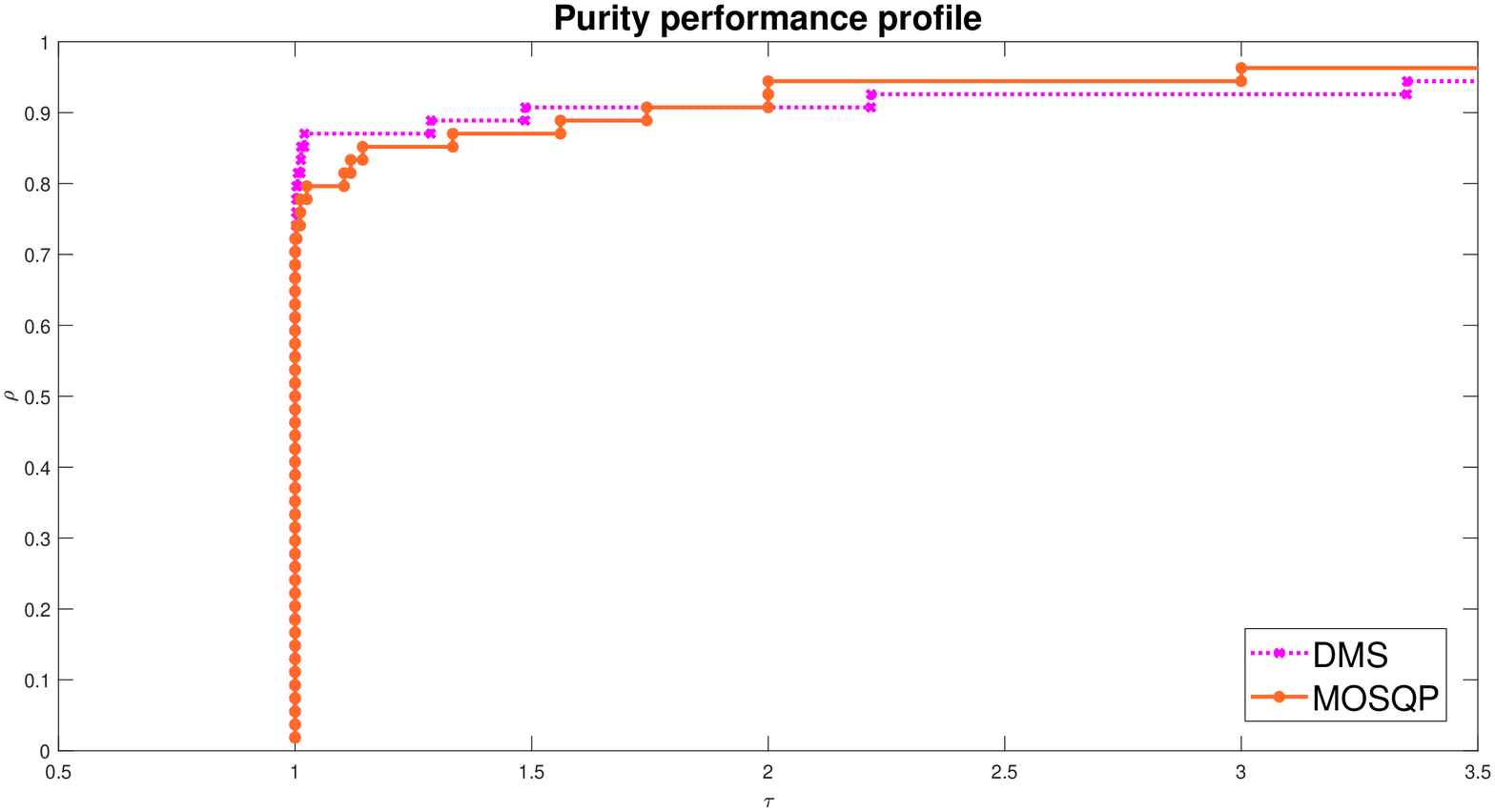}
\includegraphics[width=6.4cm,height=4cm]{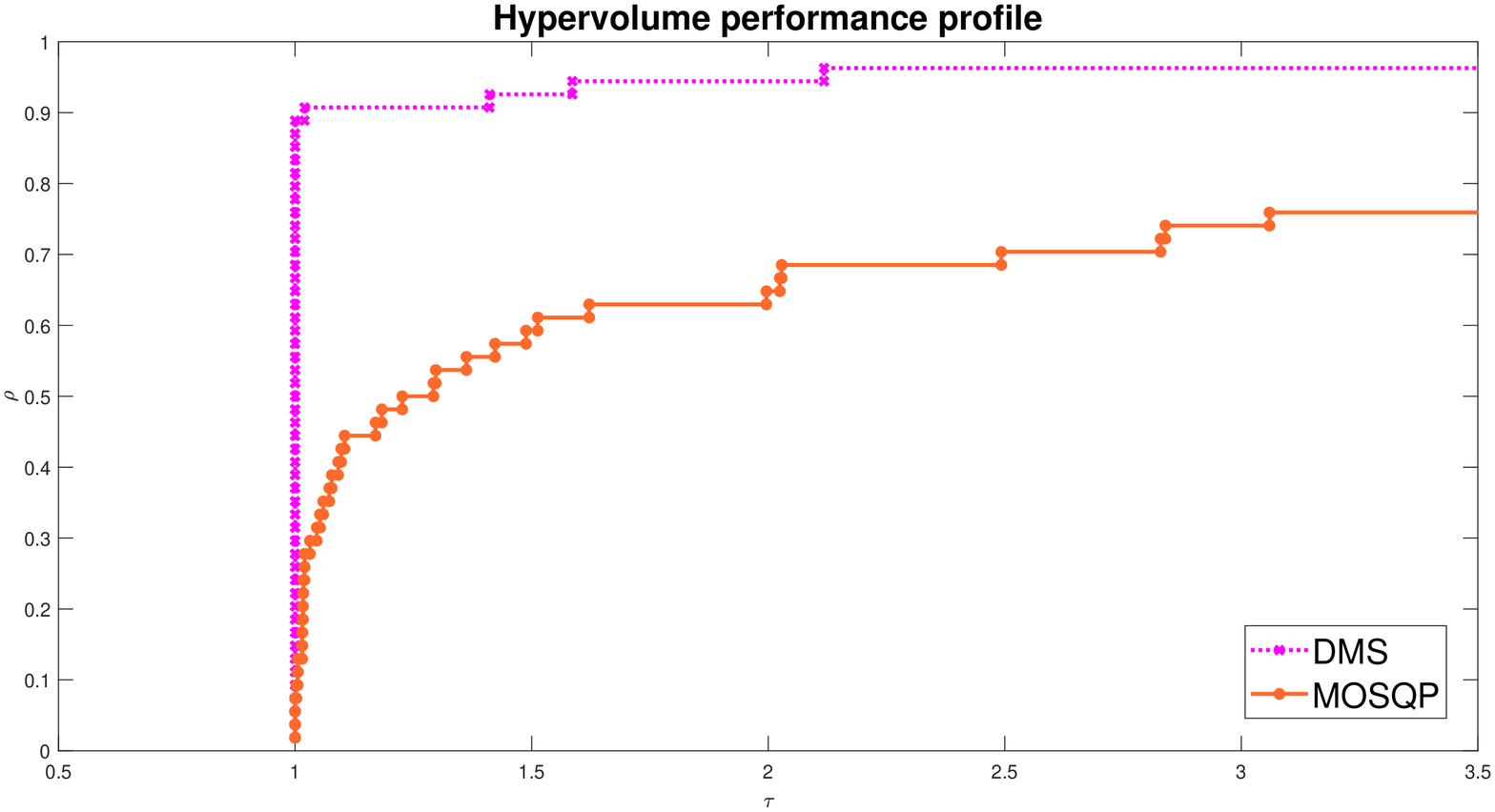}
\end{center}
\caption{\label{fig_purityhyper_DMS_vs_MOSQP}Performance profiles
for purity and hypervolume metrics, comparing the original DMS
implementation and MOSQP (maximum budget of $20\,000$ function
evaluations).}
\end{figure}

\begin{figure}[htbp]
\begin{center}
\includegraphics[width=6.4cm,height=4cm]{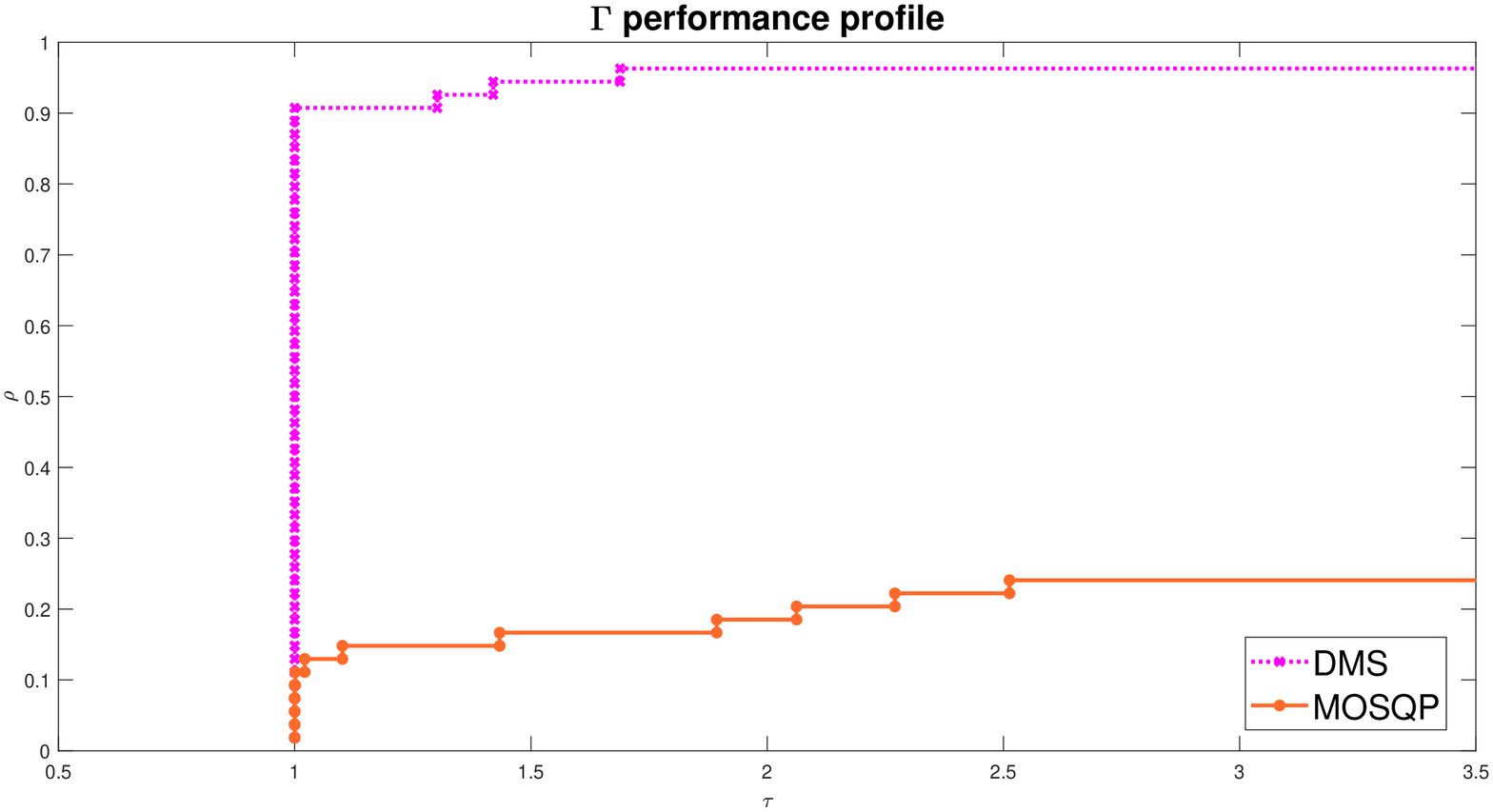}
\includegraphics[width=6.4cm,height=4cm]{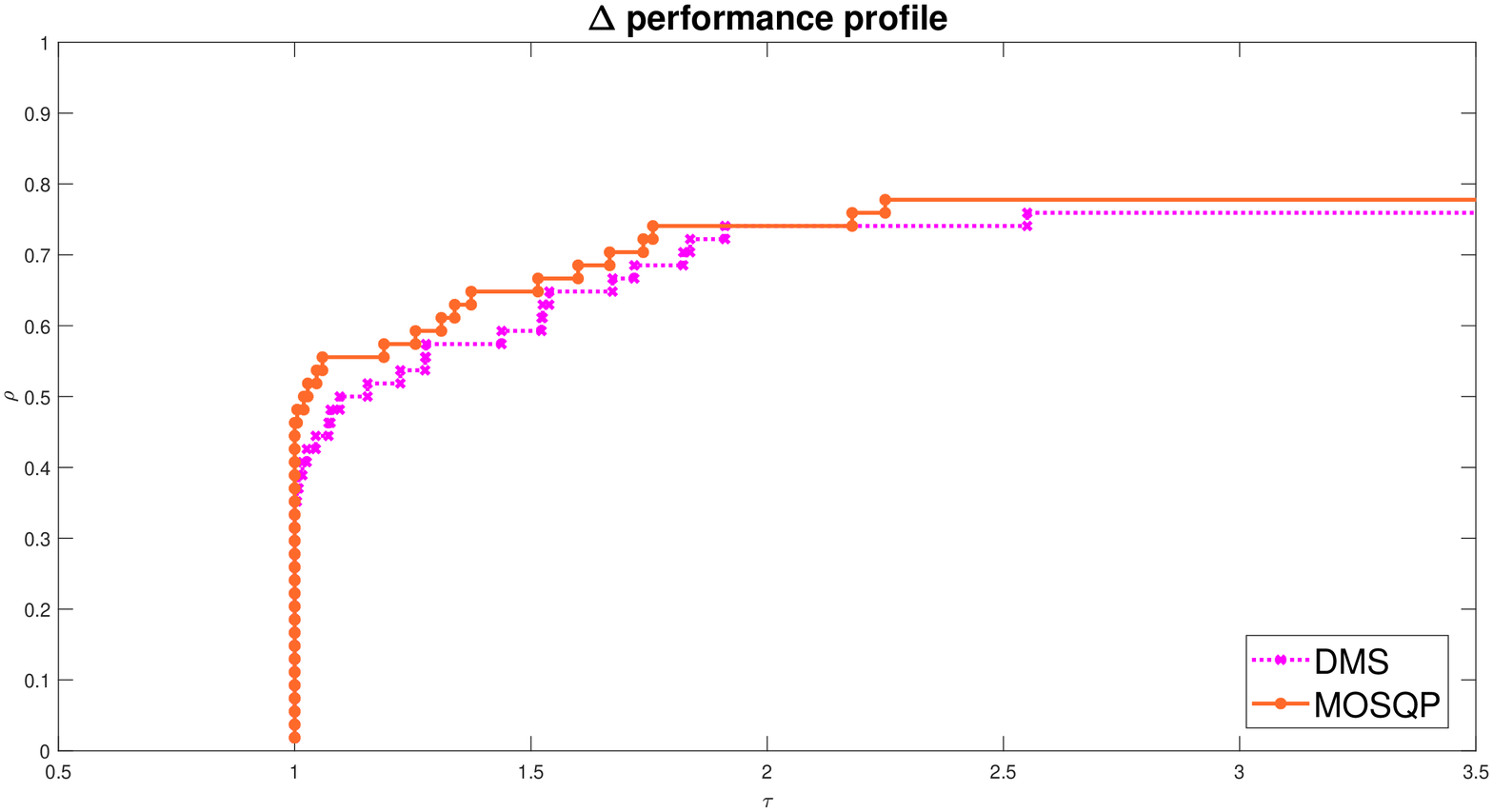}
\end{center}
\caption{\label{fig_spread_DMS_vs_MOSQP}Performance profiles for
$\Gamma$ and $\Delta$ metrics, comparing the original DMS
implementation and MOSQP (maximum budget of $20\,000$ function
evaluations).}
\end{figure}

Figures~\ref{fig_purityhyper_DMS_vs_MOSQP}
and~\ref{fig_spread_DMS_vs_MOSQP} compare DMS against MOSQP, when
a maximum budget of $20\,000$ function evaluations is considered.
The two solvers present a similar performance in terms of
robustness for purity and $\Delta$ metrics, with DMS being more
efficient in terms of purity and MOSQP with respect to the
$\Delta$ metric. However, there is a huge gain in performance with
DMS when hypervolume or the $\Gamma$ metrics are considered.

The benefits of an enriched set of directions are clear (and will
be even clearer in the following sections). If derivative-based
solvers are preferable to derivative-free optimization methods for
single objective optimization, the same does not necessarily apply
to multiobjective optimization.

\section{Pruning the poll set} \label{Sec:Pruning}
At each iteration of DMS a positive spanning set is selected as
poll set. The poll points correspond to directions in the poll
step scaled by the stepsize parameter. The objective function will
then be evaluated at all the feasible poll points, independently
of corresponding or not to descent directions.

The following result is well-known for positive spanning sets (see
Theorem 2.3 in~\cite{ARConn_et_al_2009}).
\begin{theorem}\label{theo:pss}
If $\{v_1,\ldots,v_r\}$, with $v_j\neq 0$ for all
$j\in\{1,\ldots,r\}$ positively spans $\mathbb{R}^n$ then for
every vector $d\in\mathbb{R}^n$ there is an index
$j\in\{1,\ldots,r\}$ such that $d^{\top}v_j>0$.
\end{theorem}

Considering strict differentiability of each component of the
objective function $F$, and setting $d=\nabla f_i(x)$ or
$d=-\nabla f_i(x)$, for $i\in\{1,\ldots,p\}$,
Theorem~\ref{theo:pss} allows us to conclude that in every
positive spanning set, for each component of the objective
function, we can find at least one ascent and one descent
direction.

Thus, at each iteration, for $i\in\{1,\ldots,p\}$, if $\nabla
f_i(x_k)\neq 0$, $d_k=-\nabla f_i(x_k)$ can be used to prune the
positive spanning set, only keeping directions that are descent
according to at least one component of the objective function.
Since we are only discarding directions that are ascent according
to all components of the objective function, the convergence
results of Section~\ref{Sec:DMS} still hold. The pruned set of
directions, $D^P_k$, to be considered as poll directions for DMS
at Step 2 of Algorithm 1, will then be:

$$D^P_k=\bigcup_{i\in\{1,\ldots,p\}}\{d\in D_k:-\nabla f_i(x_k)^{\top}d>0\}.$$

This strategy can be ineffective when the number of components of
the objective function is high ($p>3$), what is commonly known as
many-objective optimization~\cite{SChand_MWagner_2015}. However,
for problems comprising two or three components in the objective
function, it could lead to considerable savings in function
evaluations at each poll step.

The idea of pruning positive spanning sets was already proposed in
single objective derivative-free
optimization~\cite{MAAbramson_etal_2004}. In this setting, it is
easy to see that the cardinality of the pruned set will be $1\leq
|D^P_k|\leq |D_k|-1$. The authors were even able to provide a
particular enriched positive spanning set, that always reduces to
a singleton after pruning.

If the goal is to generate an approximation to the complete Pareto
front of a given problem, we do not wish to reduce the poll
directions to a singleton, as we do not wish to use opportunistic
approaches, which would generate at most a new feasible
nondominated point at each iteration. Moreover, in multiobjective
optimization, due to the presence of conflicting objectives, we
cannot ensure the presence of a descent direction, according to
all the components of the objective
function~\cite{ALCustodio_et_al_2011}. In fact, it is possible to
build examples where the cone of descent directions, considering
all components of the objective function, can be as narrow as one
would desire.

\begin{figure}[htbp]
\begin{center}
\includegraphics[width=6.4cm,height=4cm]{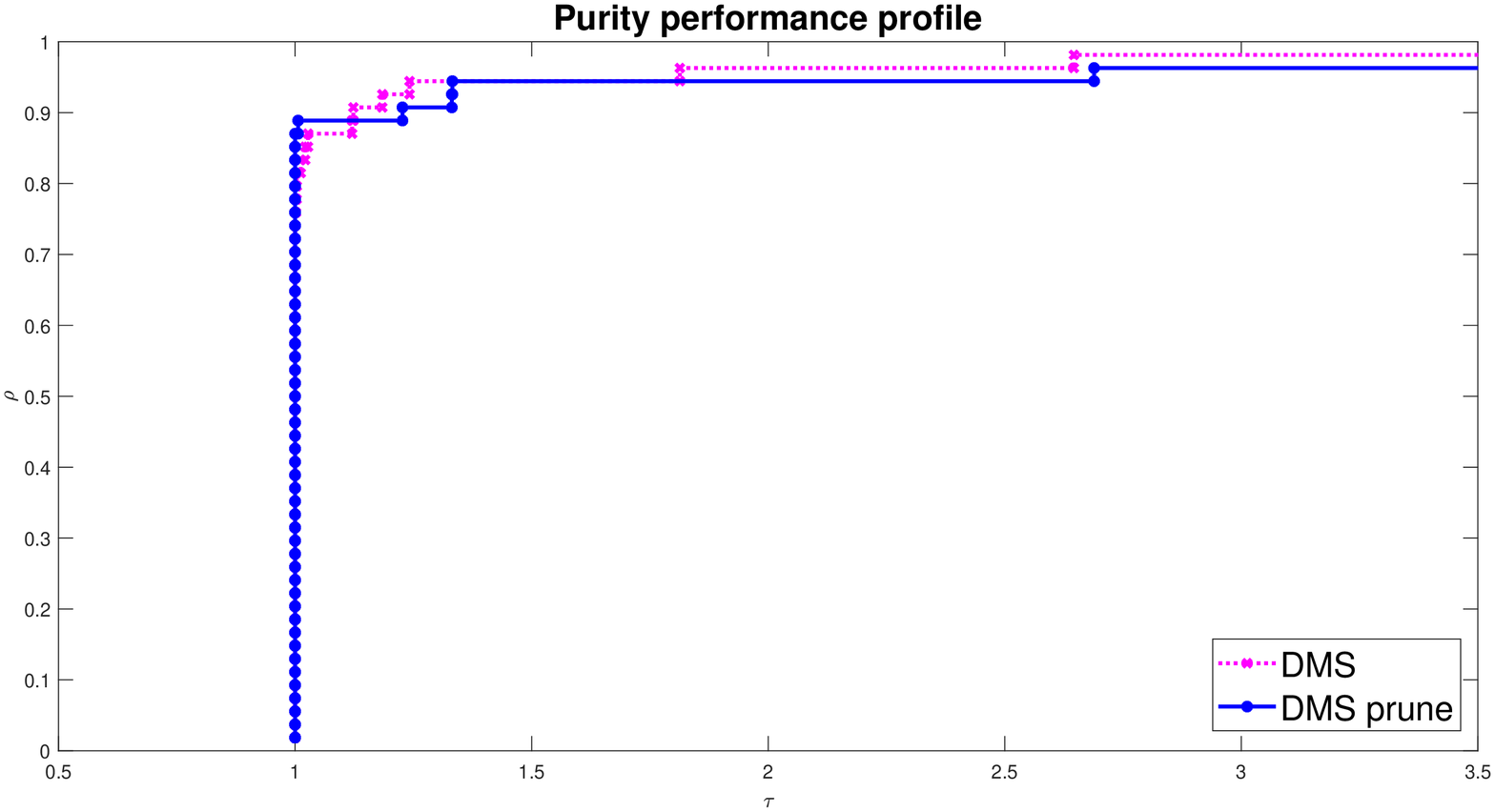}
\includegraphics[width=6.4cm,height=4cm]{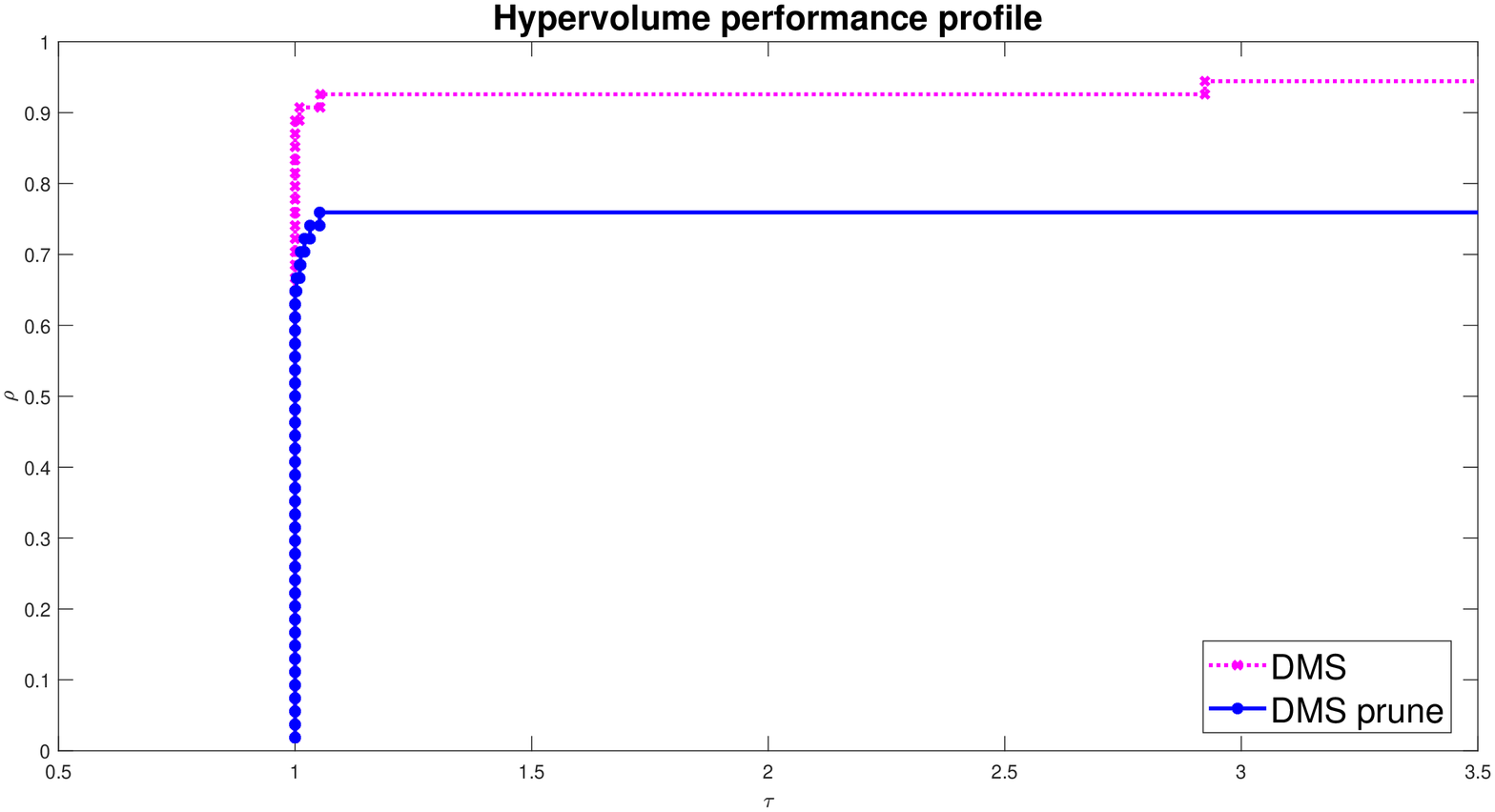}
\end{center}
\caption{\label{fig_purityhyper_DMS_vs_DMSprune}Performance
profiles for purity and hypervolume metrics, comparing the
original DMS implementation and a new version, where poll
directions are pruned using first order information (maximum
budget of $20\,000$ function evaluations).}
\end{figure}

\begin{figure}[htbp]
\begin{center}
\includegraphics[width=6.4cm,height=4cm]{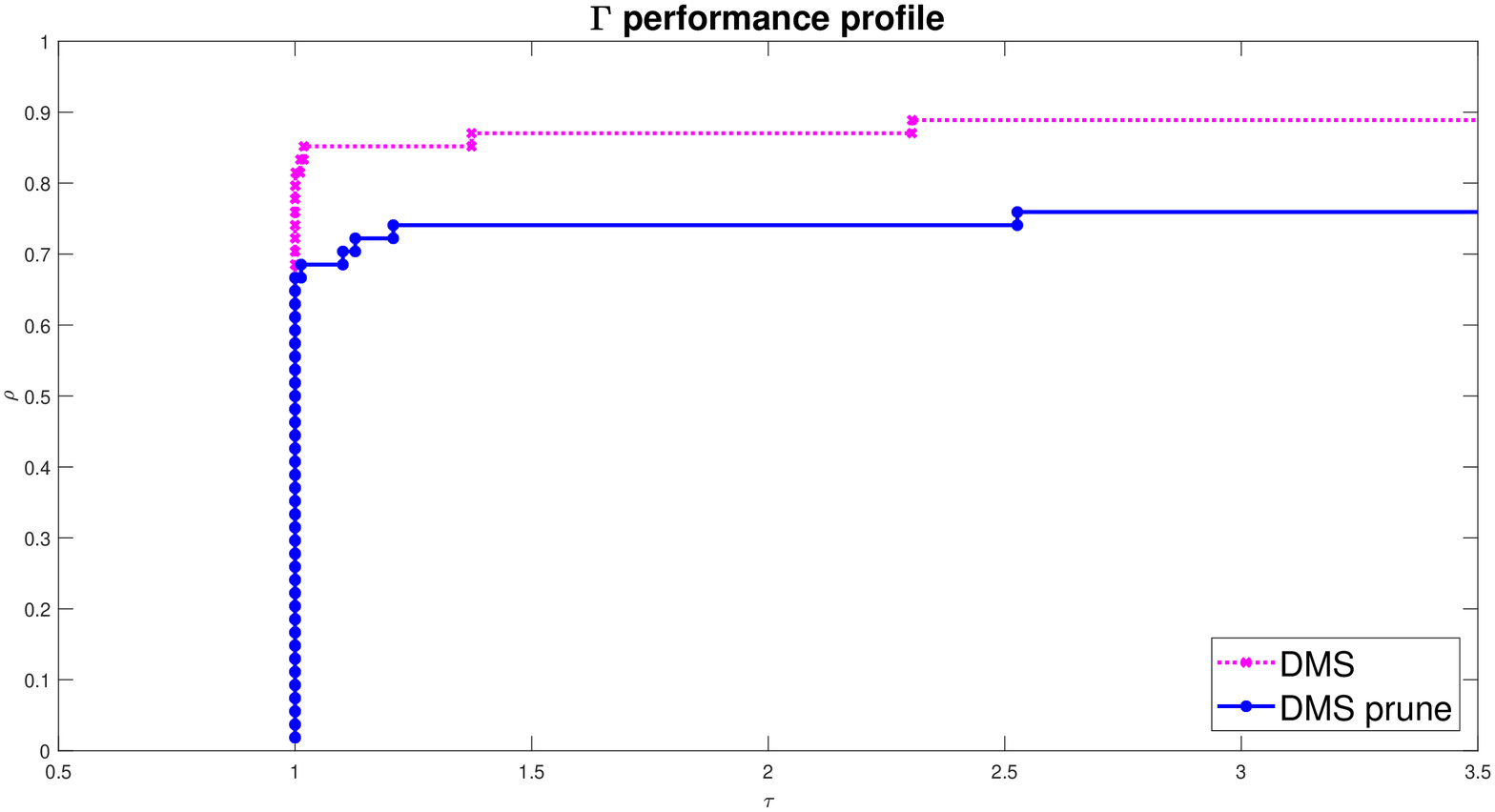}
\includegraphics[width=6.4cm,height=4cm]{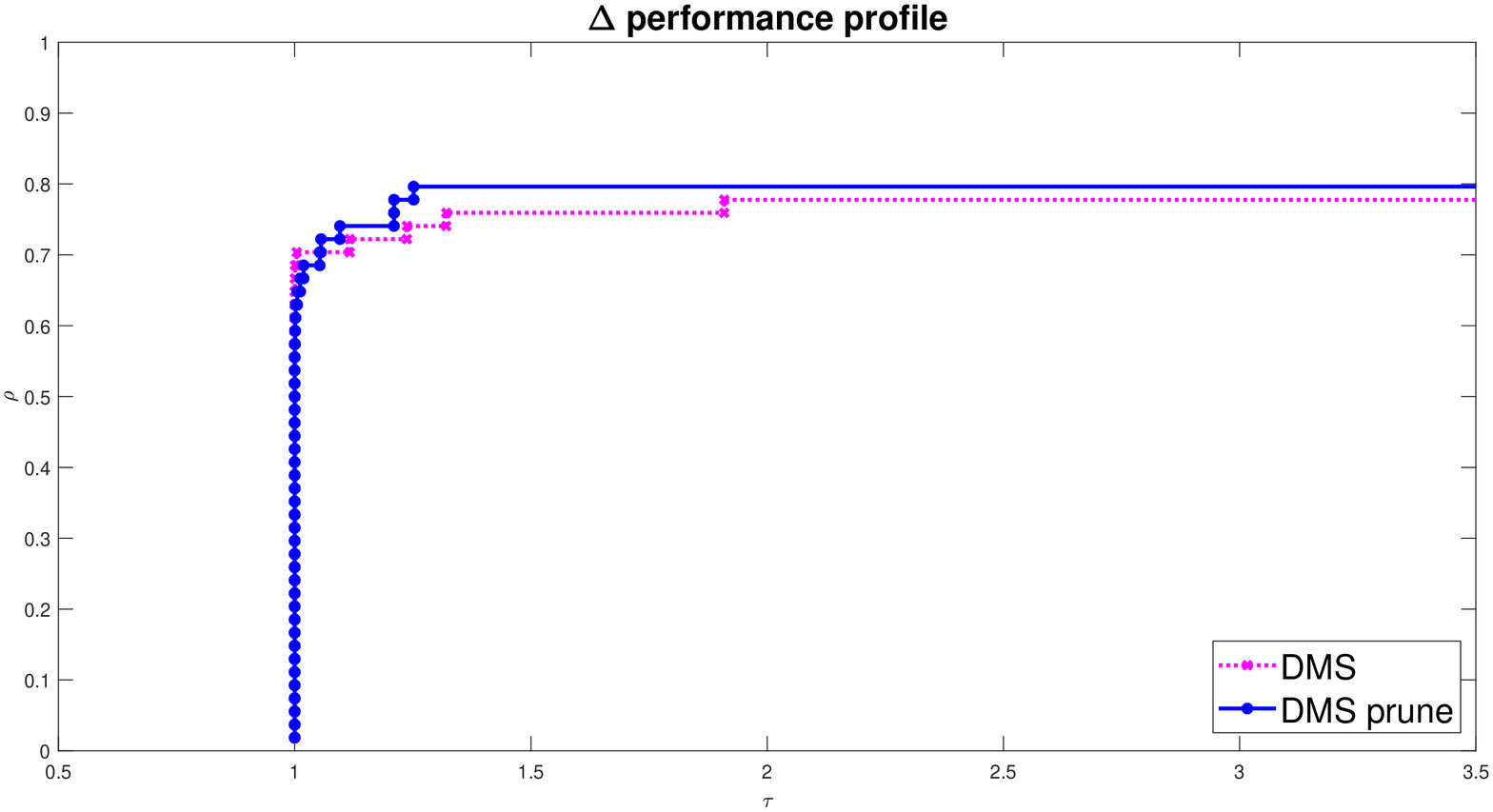}
\end{center}
\caption{\label{fig_spread_DMS_vs_DMSprune}Performance profiles
for $\Gamma$ and $\Delta$ metrics, comparing the original DMS
implementation and a new version, where poll directions are pruned
using first order information (maximum budget of $20\,000$
function evaluations).}
\end{figure}

The proposed strategy was implemented and numerically tested
against the original DMS algorithm~\cite{ALCustodio_et_al_2011}.
Figures~\ref{fig_purityhyper_DMS_vs_DMSprune}
and~\ref{fig_spread_DMS_vs_DMSprune} report the corresponding
comparison, again considering the budget of $20\,000$ function
evaluations.

In its current form, it is clear that the pruning strategy is not
successful. In fact, there is a considerable decrease in
performance regarding the hypervolume and $\Gamma$ metrics. While
it could seem surprising, as we will see in
Section~\ref{Sec:Ascent}, ascent directions play an important role
when the goal is to compute an approximation to the complete
Pareto front of a given MOO problem.

\section{The role of ascent directions} \label{Sec:Ascent}

\begin{figure}[htbp]
\begin{center}
\includegraphics[width=6.4cm,height=4cm]{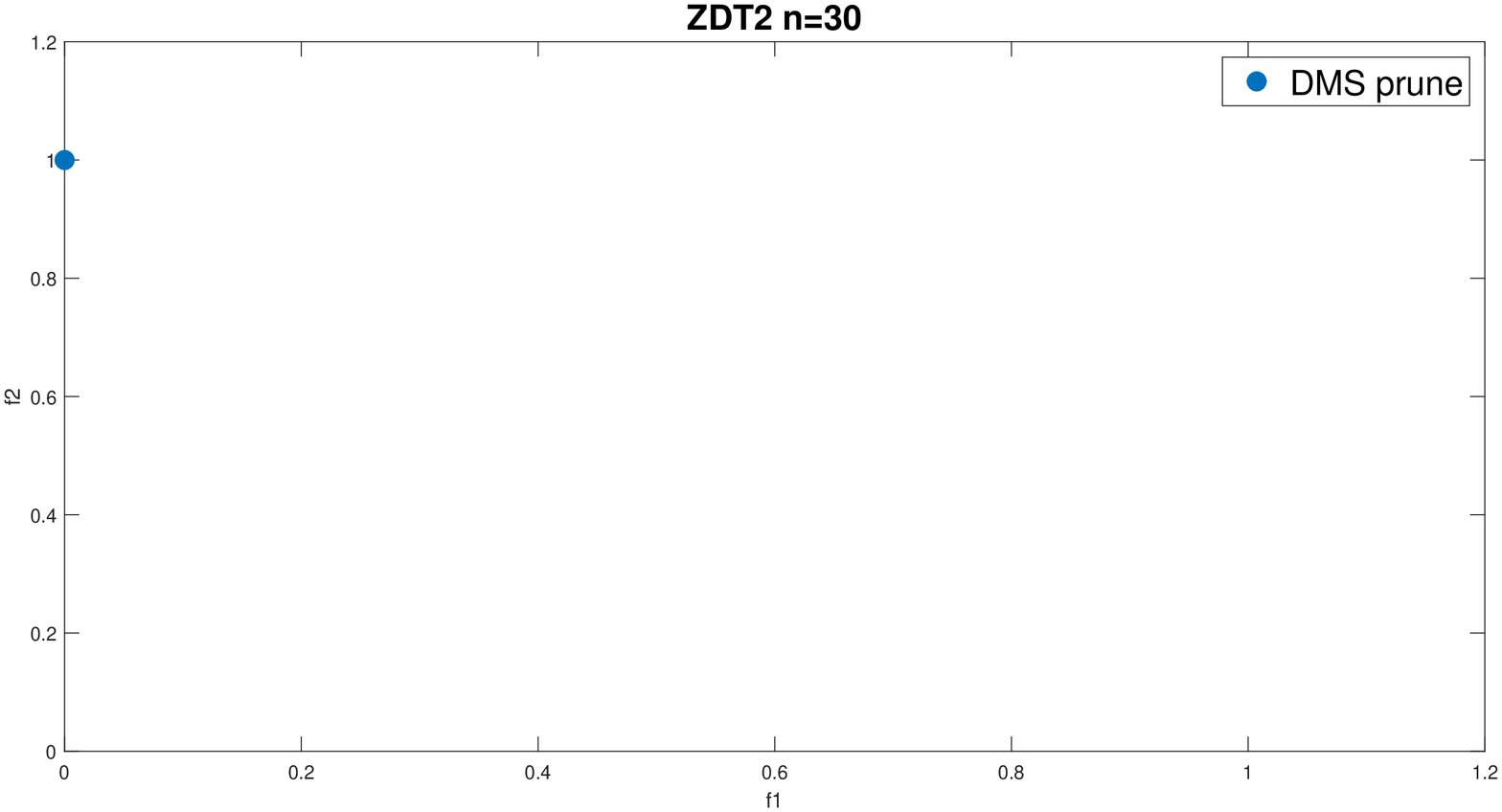}
\includegraphics[width=6.4cm,height=4cm]{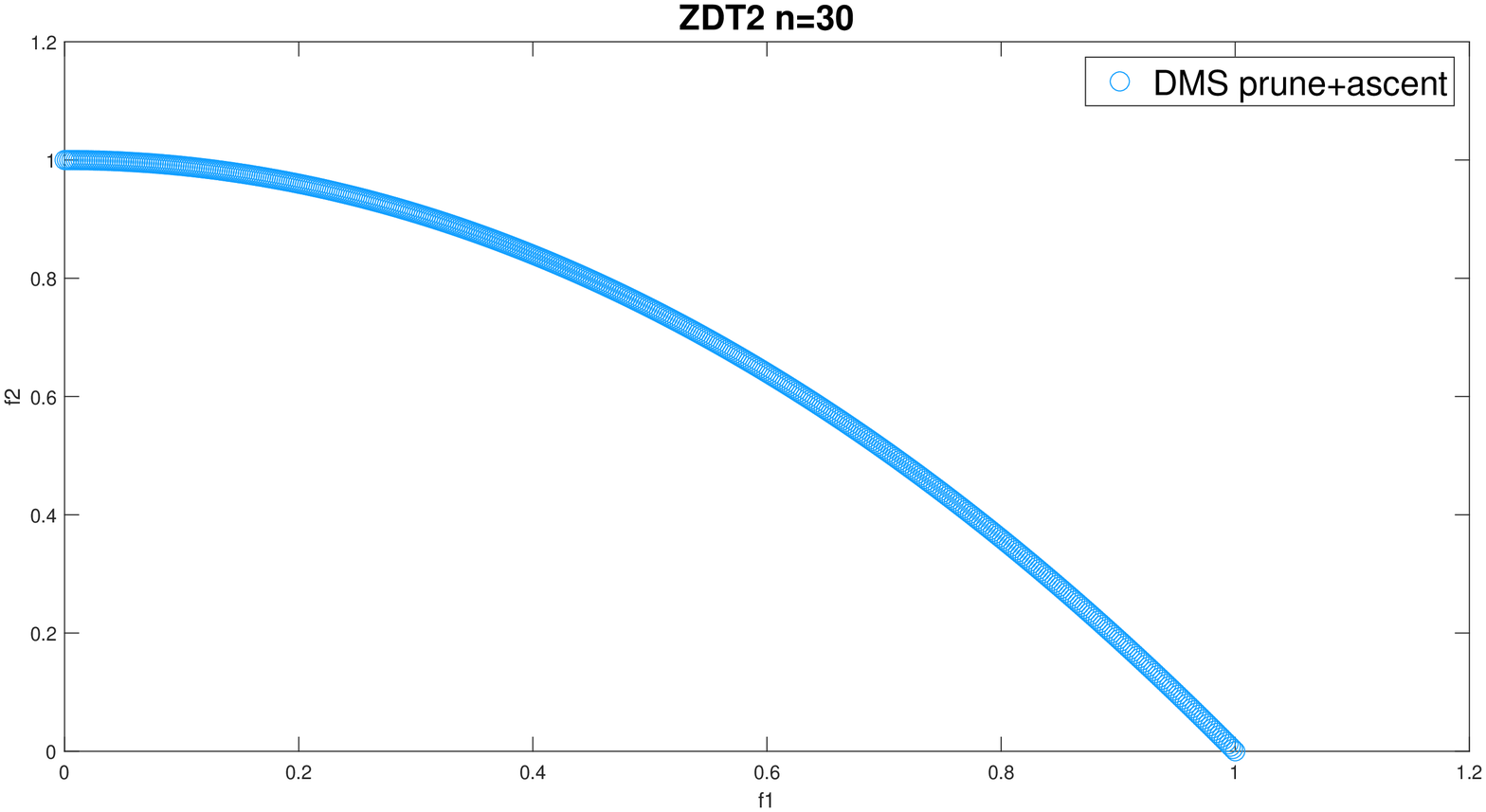}
\end{center}
\caption{\label{ZDT2}Final approximations to the Pareto front of
problem ZDT2, generated by two different algorithmic variants of
DMS. On the left, positive spanning sets are pruned to sets only
comprising descent directions. On the right, ascent directions are
considered at some iterations.}
\end{figure}

In the presence of constraints, pruning ascent directions is not
always a good strategy. Consider the biobjective minimization
problem ZDT2, with $n=30$~\cite{EZitzler_KDeb_LThiele_2000}. In
this case, if we provide as initialization one Pareto critical
point, the algorithm stalls, being unable to generate other Pareto
critical points in the Pareto front. This behavior is accordingly
to the convergence results derived for DMS, which only guarantee
convergence to a single Pareto critical point. By providing ascent
directions, that conform to the geometry of the nearby feasible
region, the algorithm is able to proceed and generate a large
number of Pareto critical points. Figure~\ref{ZDT2} illustrates
the situation.

\begin{figure}[htbp]
\begin{center}
\includegraphics[width=6.4cm,height=4cm]{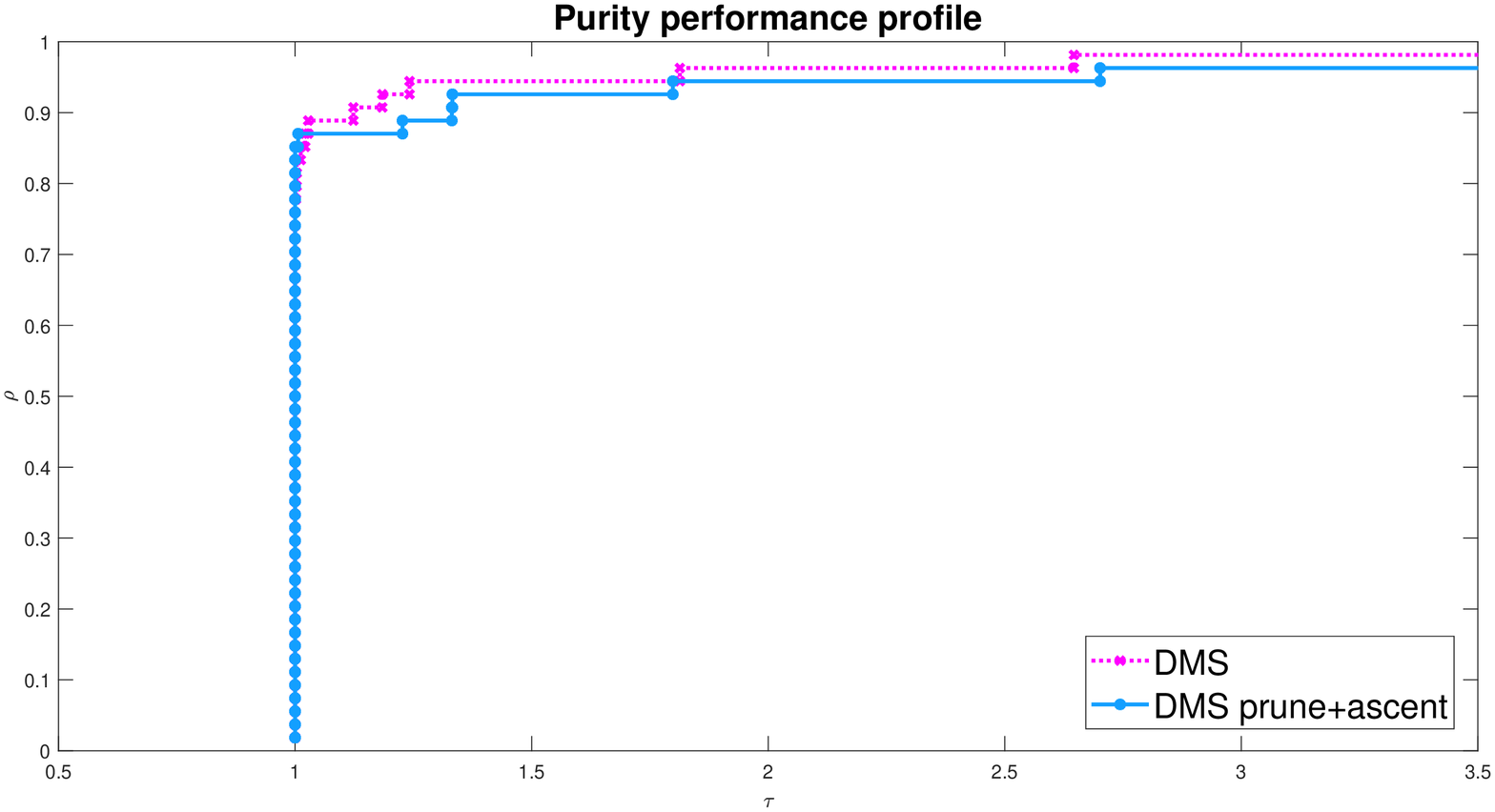}
\includegraphics[width=6.4cm,height=4cm]{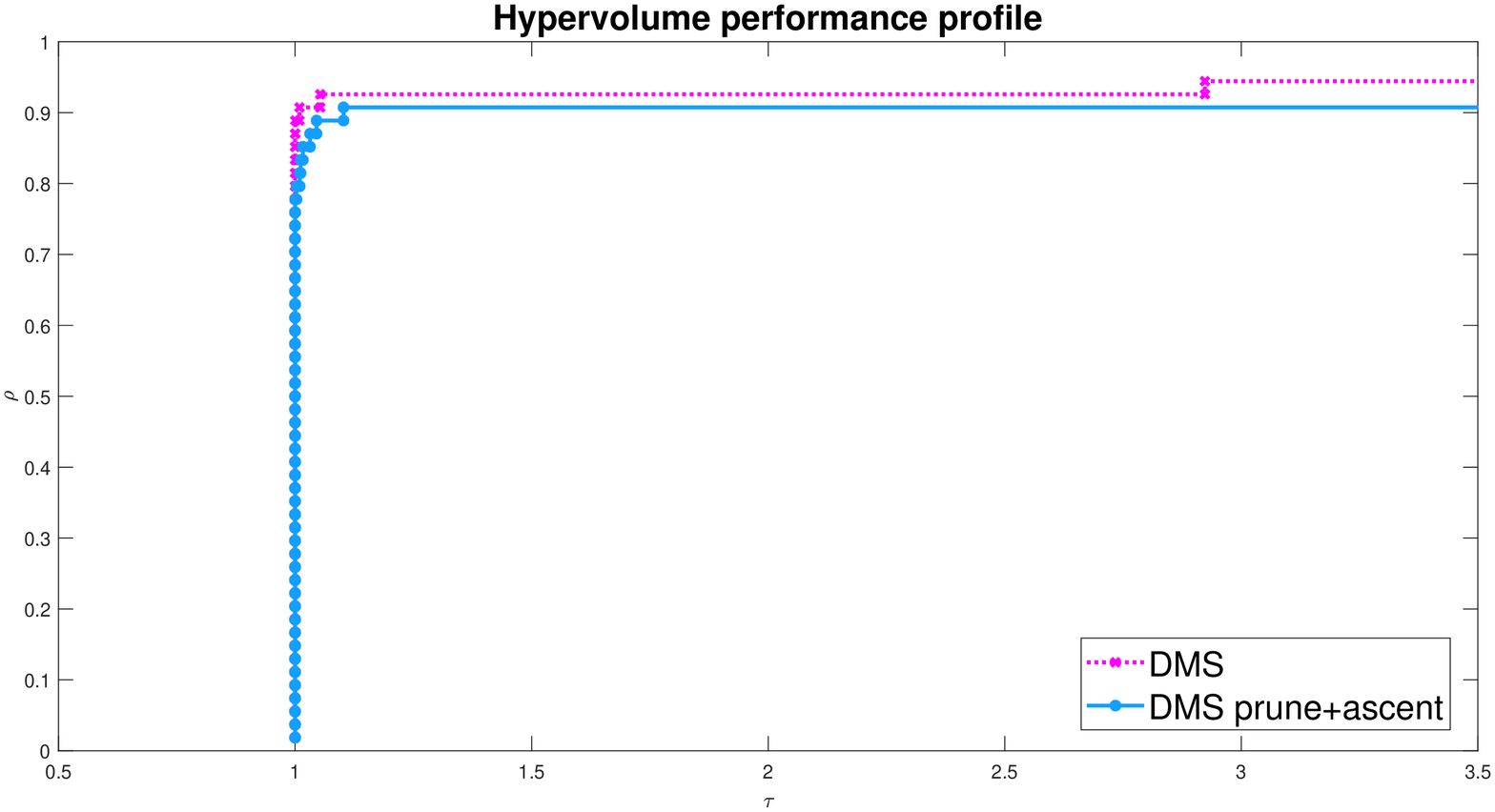}
\end{center}
\caption{\label{fig_purityhyper_DMS_vs_DMSpruneascent}Performance
profiles for purity and hypervolume metrics, comparing the
original DMS implementation and a new version, where poll
directions are pruned using first order information, but not at
all the iterations (maximum budget of $20\,000$ function
evaluations).}
\end{figure}

\begin{figure}[htbp]
\begin{center}
\includegraphics[width=6.4cm,height=4cm]{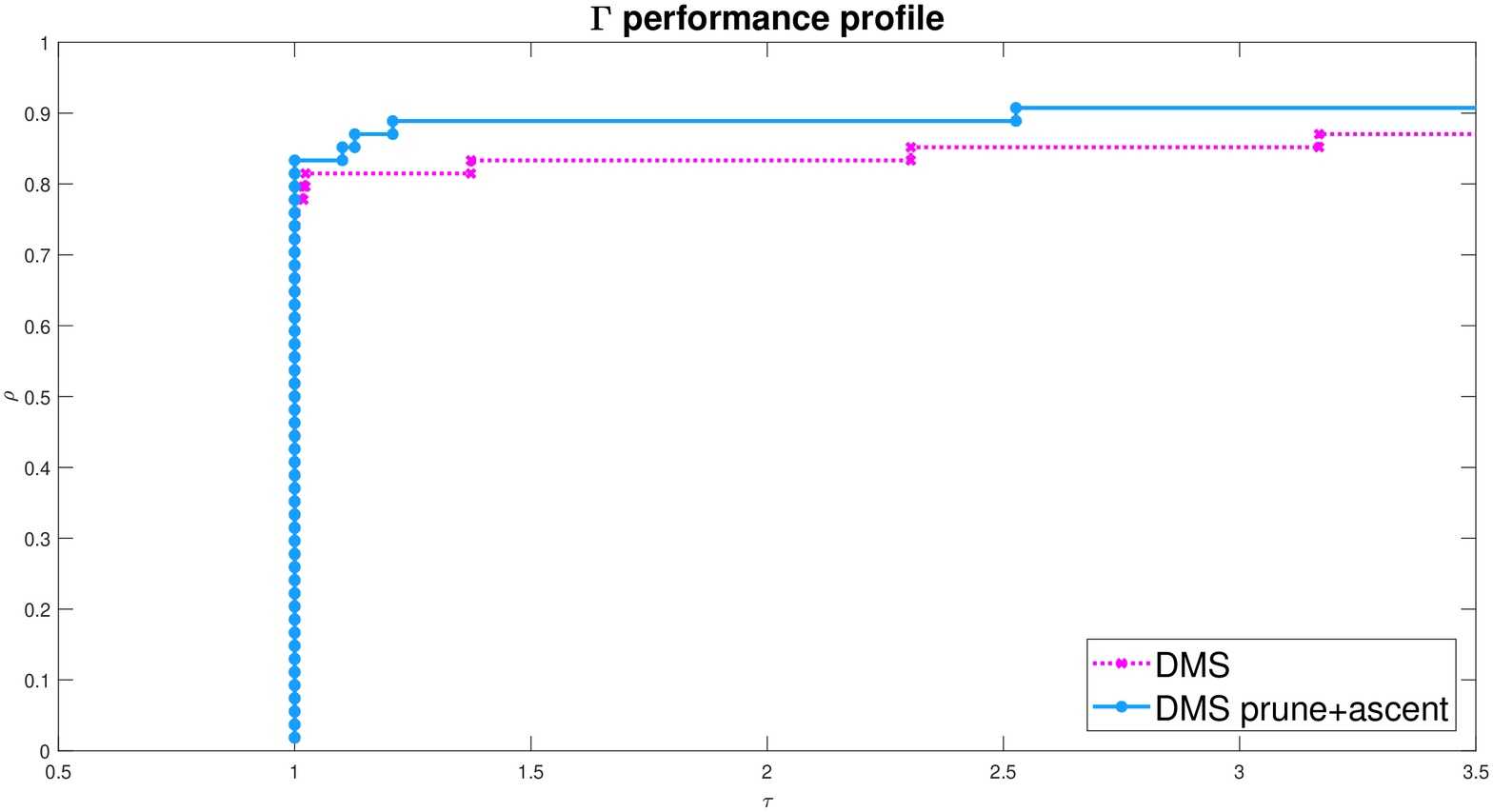}
\includegraphics[width=6.4cm,height=4cm]{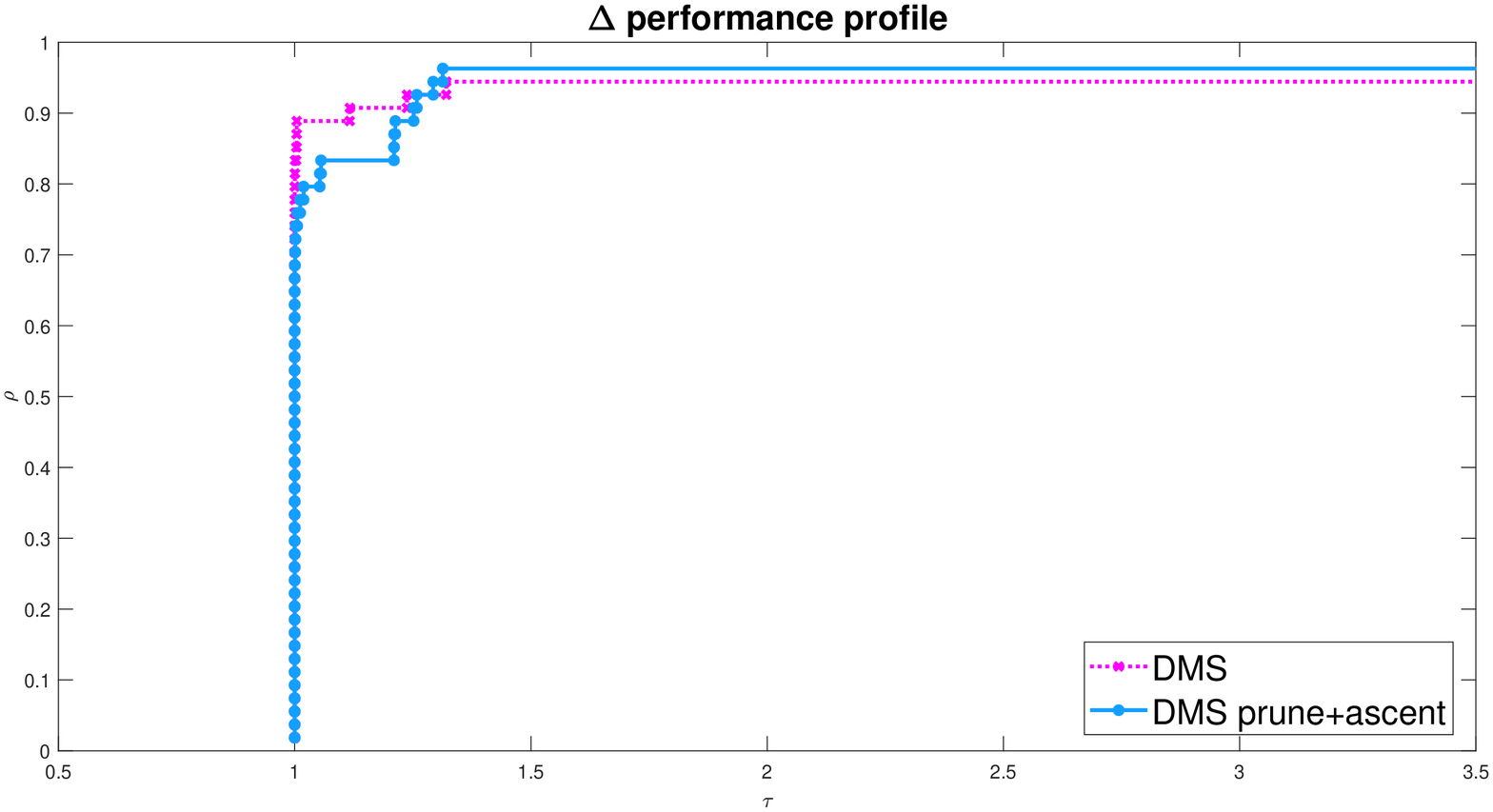}
\end{center}
\caption{\label{fig_spread_DMS_vs_DMSpruneascent}Performance
profiles for $\Gamma$ and $\Delta$ metrics, comparing the original
DMS implementation and a new version, where poll directions are
pruned using first order information, but not at all the
iterations (maximum budget of $20\,000$ function evaluations).}
\end{figure}

Thus, the approach taken was to return to the original positive
spanning set $D_k$ (without pruning) at some iterations. Assume
that at a given iteration the original positive spanning set was
pruned and $D_k^P$ was used as poll set, but the algorithm was
unable to proceed because every poll point was infeasible. At the
next iteration pruning will not be applied, and the original
positive spanning set $D_k$ will be considered as the set of poll
directions. Again, since we are only disregarding directions that
are ascent according to all components of the objective function,
and only at some iterations, the convergence results of
Section~\ref{Sec:DMS} continue to hold.
Figures~\ref{fig_purityhyper_DMS_vs_DMSpruneascent} and
\ref{fig_spread_DMS_vs_DMSpruneascent} report performance profiles
comparing this new approach with the original implementation of
DMS.

\begin{figure}[htbp]
\begin{center}
\includegraphics[width=6.4cm,height=4cm]{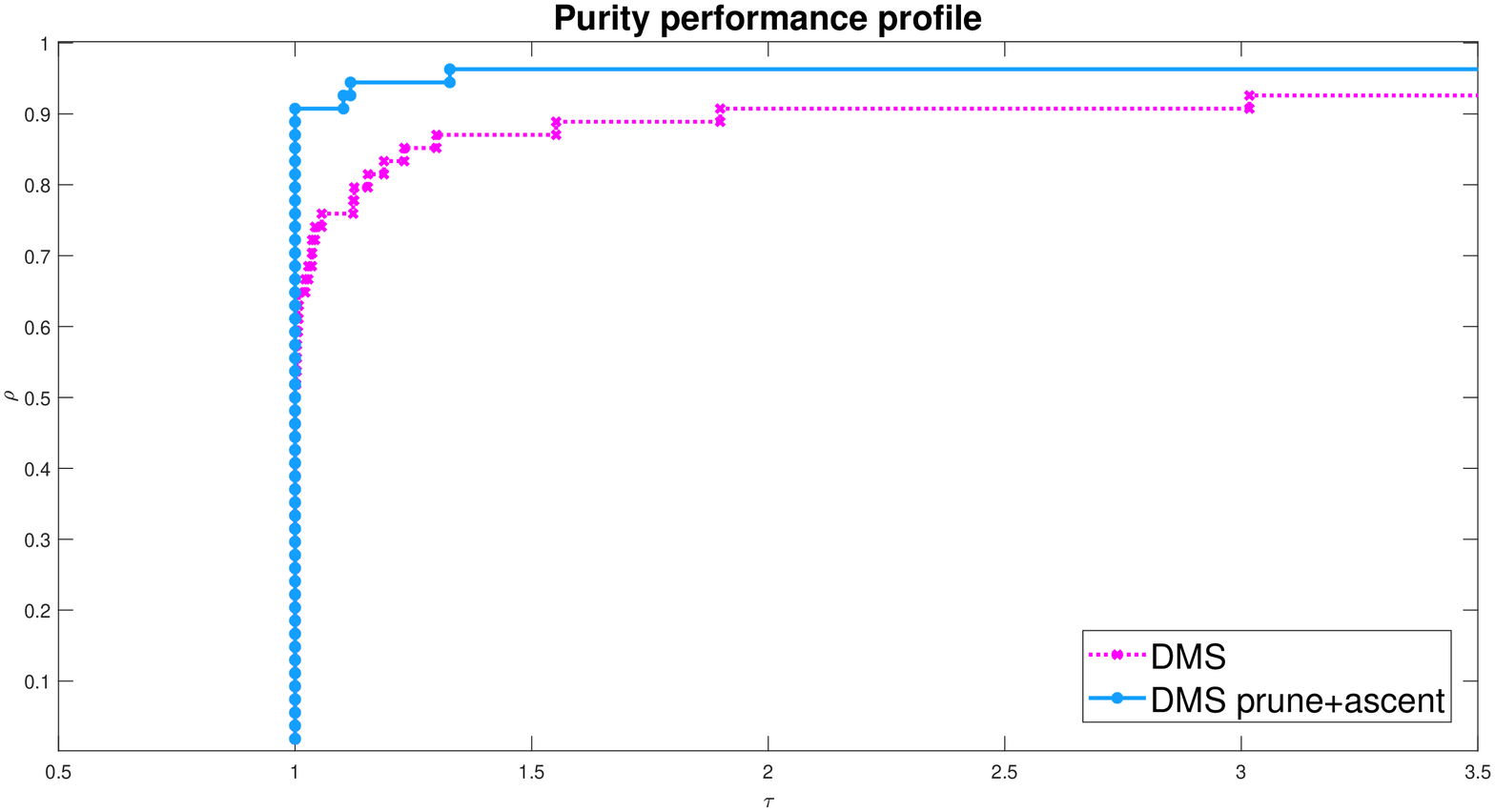}
\includegraphics[width=6.4cm,height=4cm]{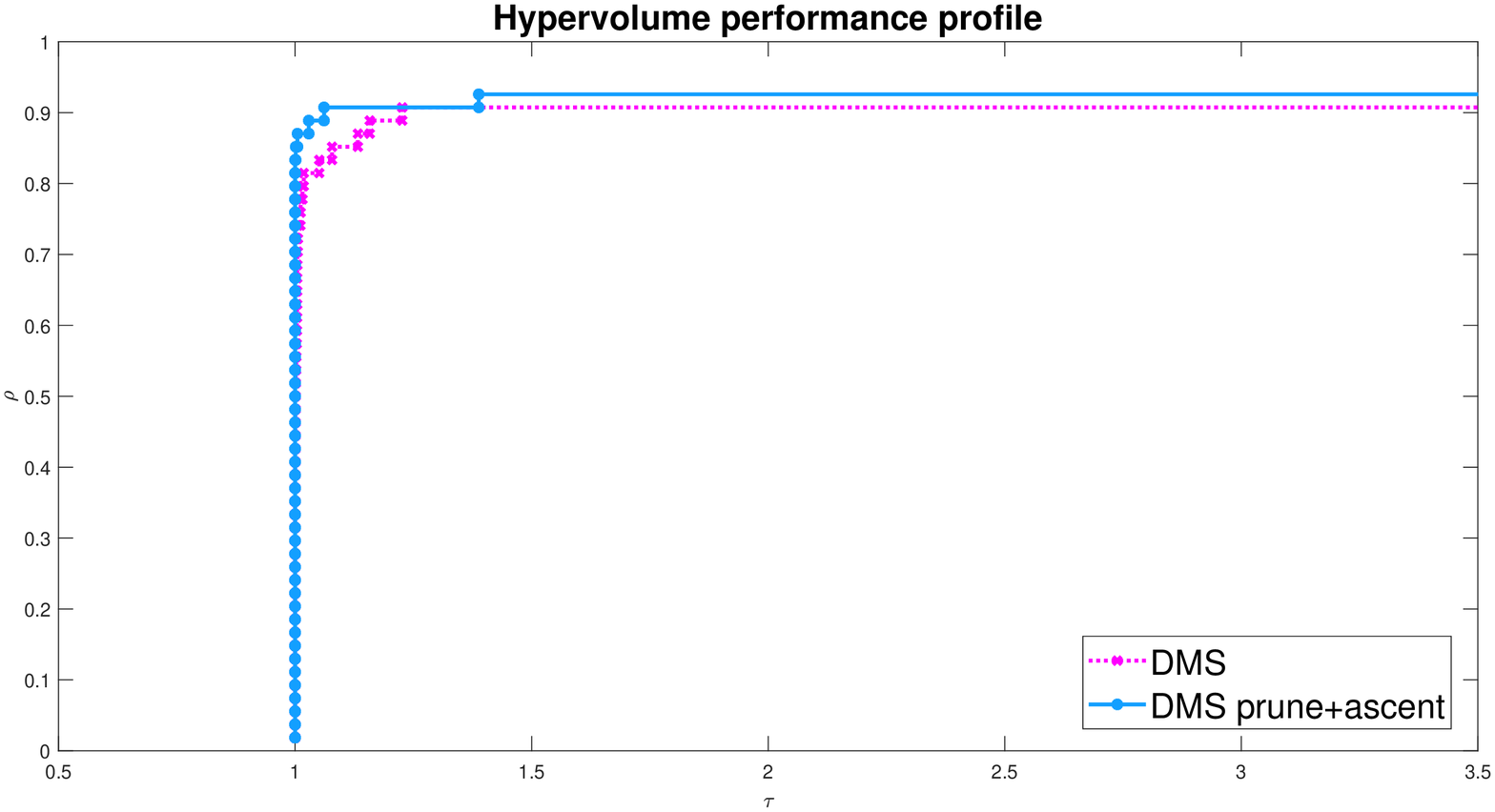}
\end{center}
\caption{\label{fig_purityhyper_DMS_vs_DMSpruneascent500}Performance
profiles for purity and hypervolume metrics, comparing the
original DMS implementation and a new version, where poll
directions are pruned using first order information, but not at
all the iterations (maximum budget of $500$ function
evaluations).}
\end{figure}

\begin{figure}[htbp]
\begin{center}
\includegraphics[width=6.4cm,height=4cm]{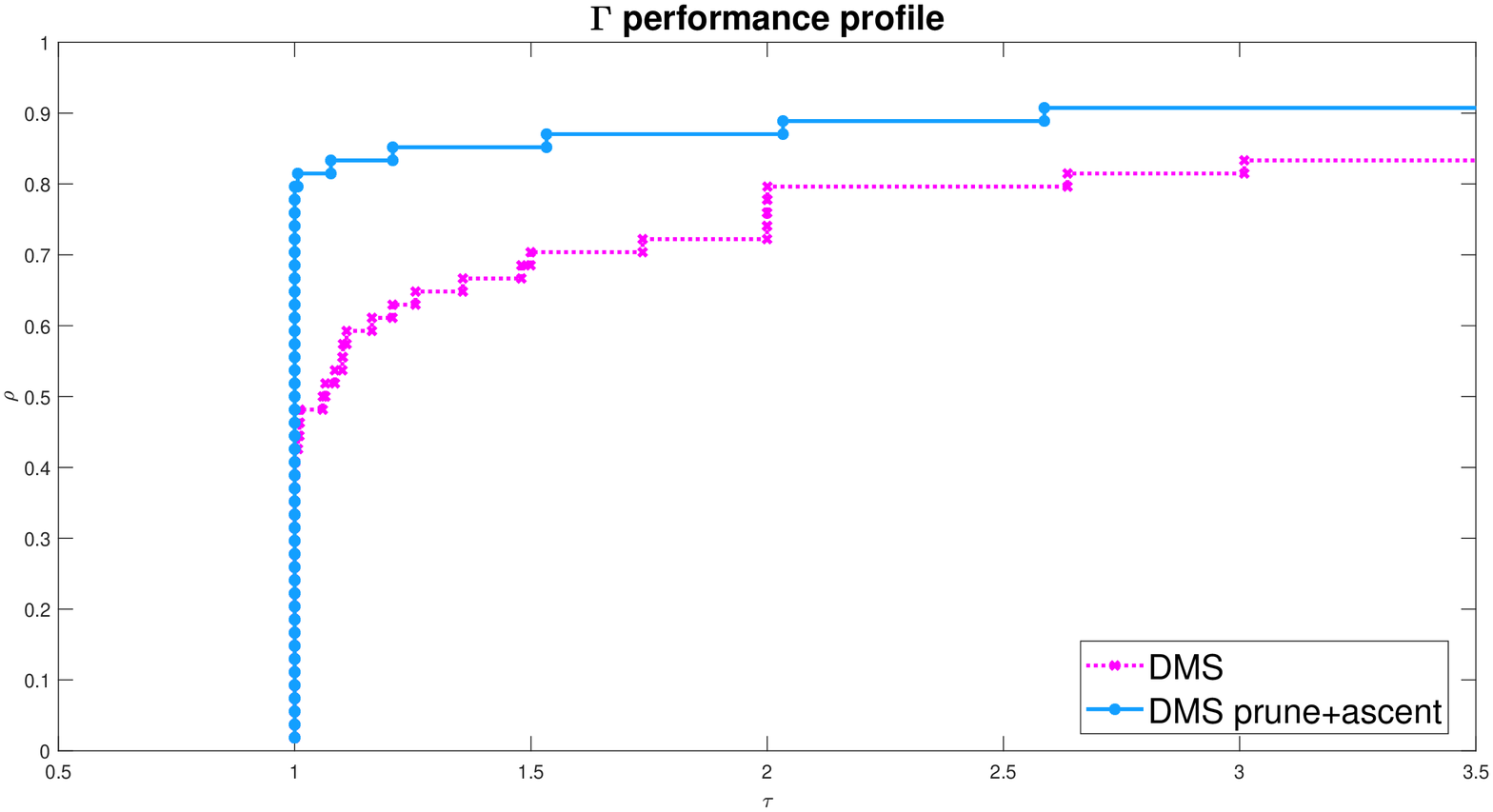}
\includegraphics[width=6.4cm,height=4cm]{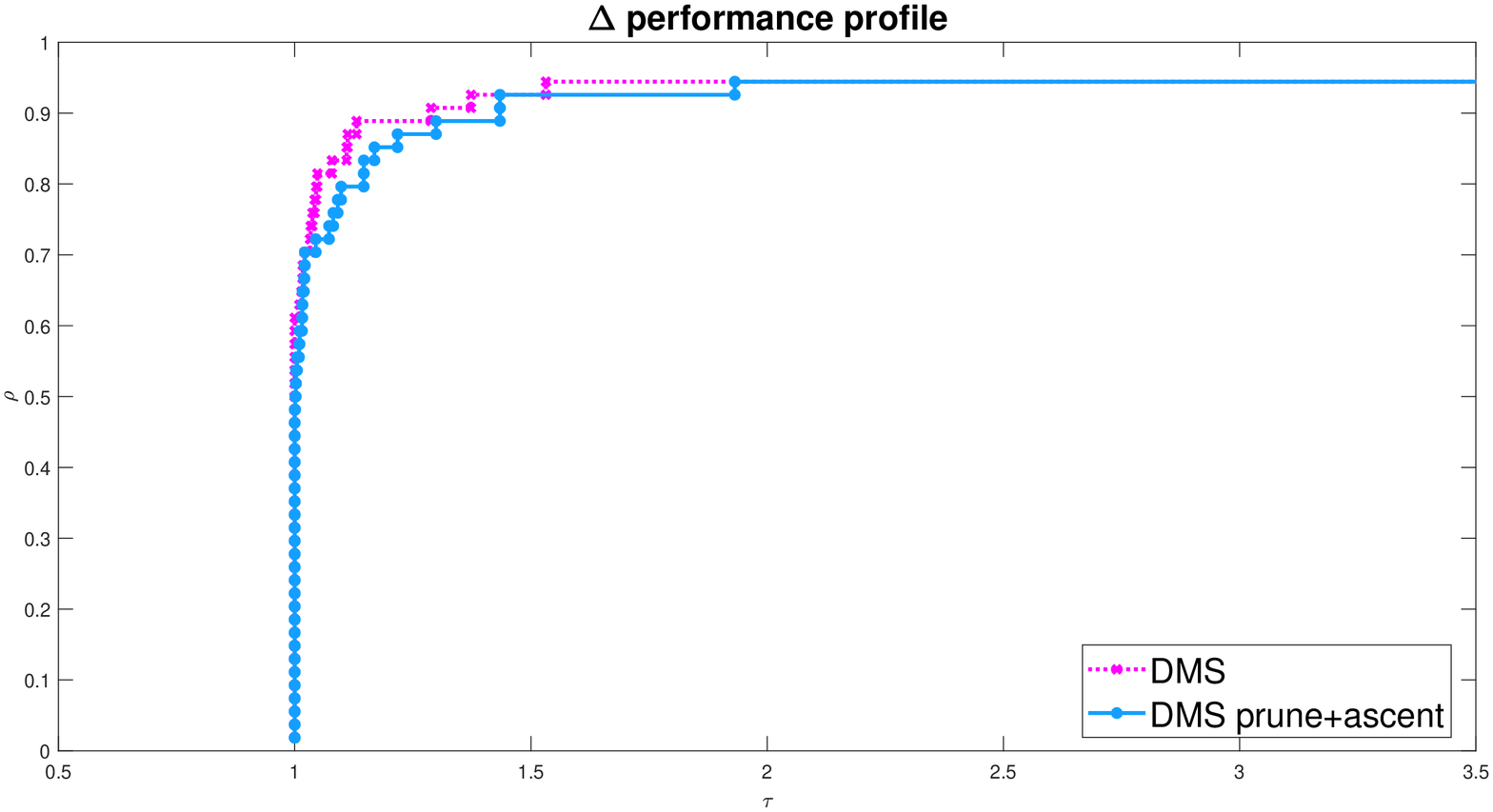}
\end{center}
\caption{\label{fig_spread_DMS_vs_DMSpruneascent500}Performance
profiles for $\Gamma$ and $\Delta$ metrics, comparing the original
DMS implementation and a new version, where poll directions are
pruned using first order information, but not at all the
iterations (maximum budget of $500$ function evaluations).}
\end{figure}

Now, the two variants of DMS are extremely close in terms of
performance, but the new approach brings some advantage in terms
of efficiency for the $\Gamma$ metric. However, the advantages of
the new approach are clearer if the computational budget is
reduced from $20\,000$ to only $500$ functions evaluations (see
Figures~\ref{fig_purityhyper_DMS_vs_DMSpruneascent500}
and~\ref{fig_spread_DMS_vs_DMSpruneascent500}).

Savings in function evaluations, allow clear improvements in terms
of efficiency, both for purity and hypervolume. Regarding the
$\Gamma$ metric, there is a clear advantage of the new variant
over the classical DMS approach.

\begin{figure}[htbp]
\begin{center}
\includegraphics[width=6.4cm,height=4cm]{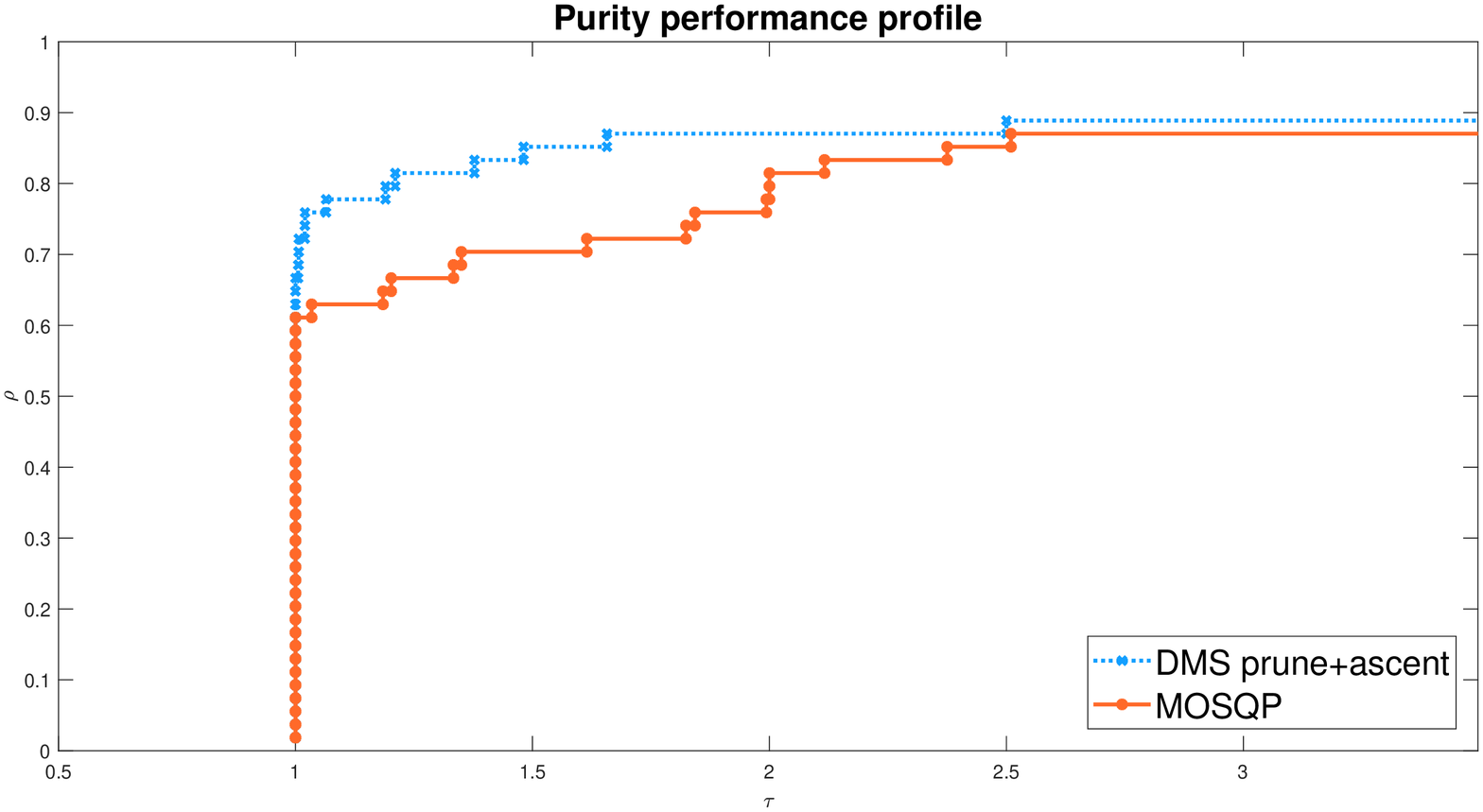}
\includegraphics[width=6.4cm,height=4cm]{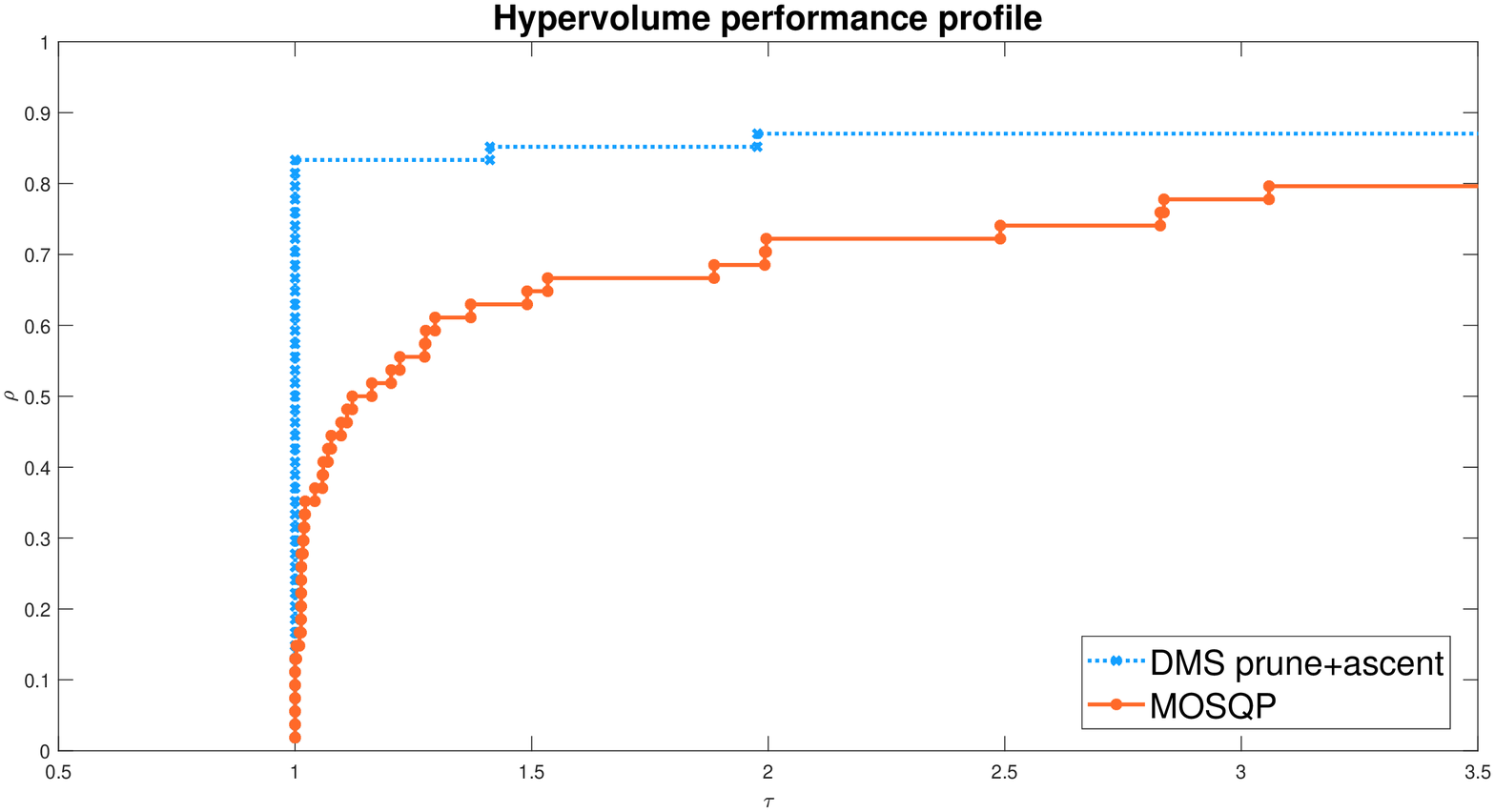}
\end{center}
\caption{\label{fig_purityhyper_MOSQP_vs_DMSpruneascent500}Performance
profiles for purity and hypervolume metrics, comparing MOSQP and
the new version of DMS, where poll directions are pruned using
first order information, but not at all the iterations (maximum
budget of $500$ function evaluations).}
\end{figure}

\begin{figure}[htbp]
\begin{center}
\includegraphics[width=6.4cm,height=4cm]{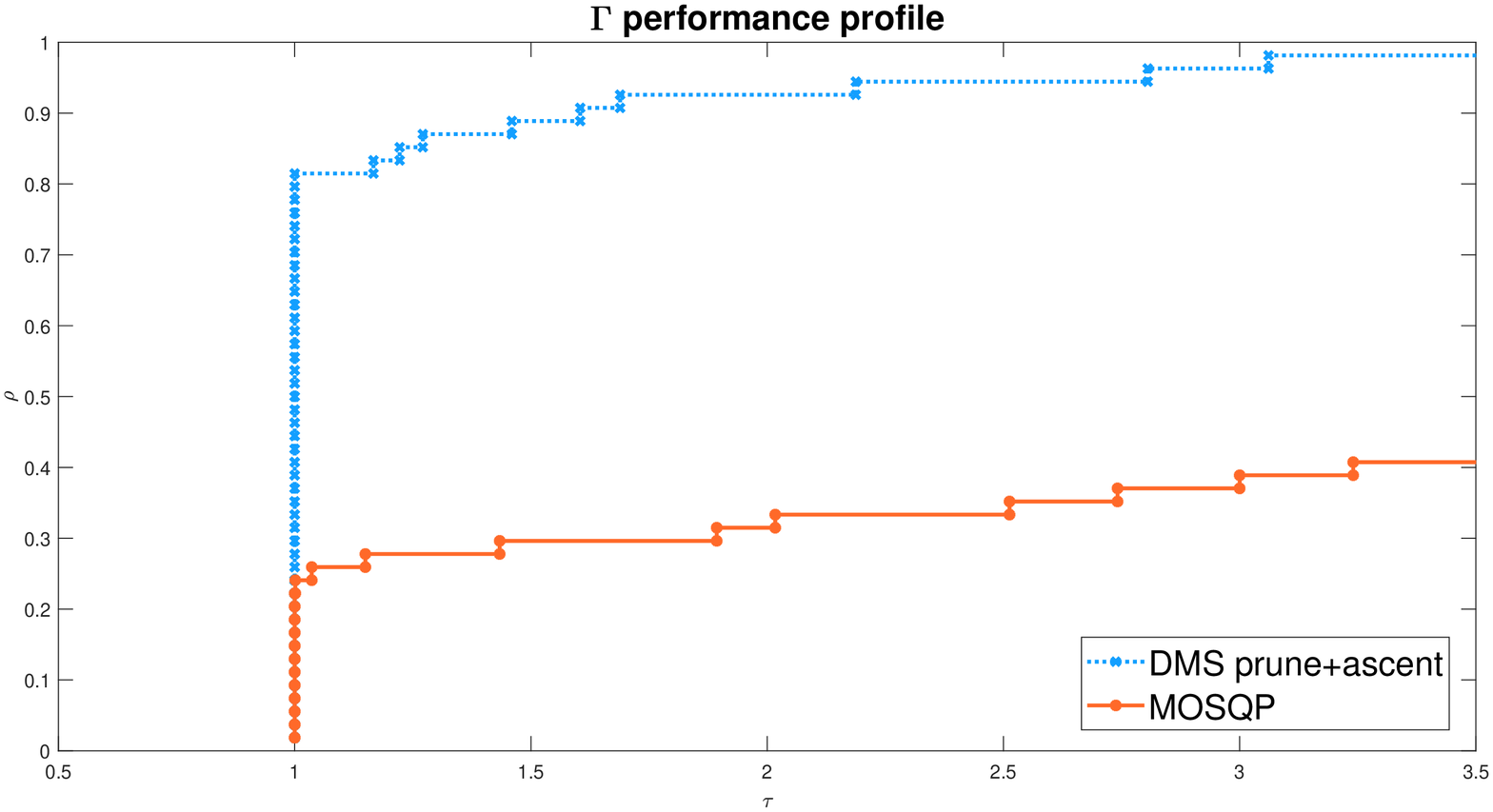}
\includegraphics[width=6.4cm,height=4cm]{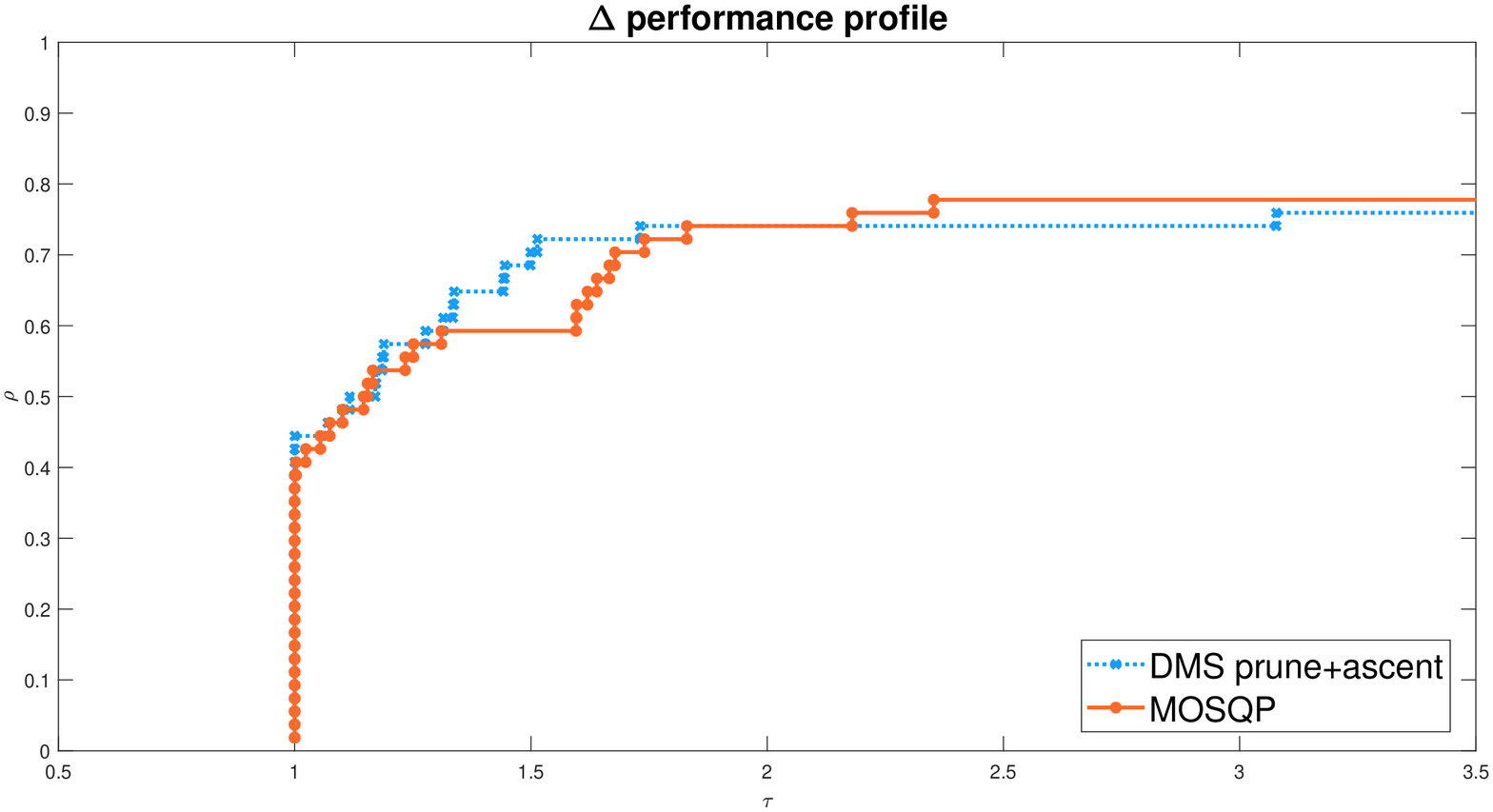}
\end{center}
\caption{\label{fig_spread_MOSQP_vs_DMSpruneascent500}Performance
profiles for $\Gamma$ and $\Delta$ metrics, comparing MOSQP and
the new version of DMS, where poll directions are pruned using
first order information, but not at all the iterations (maximum
budget of $500$ function evaluations).}
\end{figure}

Comparing with MOSQP, again considering a budget of only $500$
function evaluations, there is also a clear advantage of the new
variant in three of the metrics considered, namely purity,
hypervolume and $\Gamma$ metric, with an equal performance for the
$\Delta$ metric.
Figures~\ref{fig_purityhyper_MOSQP_vs_DMSpruneascent500}
and~\ref{fig_spread_MOSQP_vs_DMSpruneascent500} report the
results. Thus, the good performance of derivative-free
optimization solvers over derivative-based ones is not the result
of large budgets of function evaluations, but a consequence of
richer sets of directions, including ascent ones.

\section{Concluding remarks}\label{Sec:Conclusions}

DMS was proposed in~\cite{ALCustodio_et_al_2011} as a robust and
efficient algorithm, able to generate approximations to the
complete Pareto front of MOO problems. Surprisingly, it showed to
be a strong competitor against the derivative-based solver MOSQP,
evidencing that in MOO, when the goal is to generate an
approximation to the complete Pareto front of a given problem,
even if first-order derivatives are available, derivative-free
solvers can be good alternatives to derivative-based approaches.

Derivatives can be used to prune the positive spanning sets to be
considered as poll directions. However, care should be taken
because ascent directions, that conform to the geometry of the
nearby feasible region, can have an important role in the ability
of generating a complete approximation to the Pareto front of a
given problem.

The new variant of DMS, which prunes the poll set of directions,
but that at some iterations considers its enrichment with ascent
directions, showed to be competitive both with the
derivative-based solver MOSQP and with the original implementation
of DMS. For low computational budgets of function evaluations, it
allows an increase in the percentage of nondominated points
generated in the approximation to the Pareto front of the MOO
problem and also a reduction in the largest gap across the
generated Pareto front, when compared with the original
implementation of DMS. In the case of MOSQP, there are additional
advantages regarding the hypervolume associated to the computed
approximation to the Pareto front.

Future work could include the definition of a search step taking
advantage of first order information for building Taylor models,
which will be minimized considering an approach similar to the one
proposed and analyzed in~\cite{CPBras_ALCustodio_2020}.

\small

\bibliographystyle{siam}
\bibliography{ref-EnrichedDMS}

\end{document}